\documentclass{irmaart}
\pdfoutput=1

\usepackage{IRMAstyle_alex}

\newcommand*\circled[1]{\tikz[baseline=(char.base)]{
            \node[shape=circle,draw,inner sep=2pt] (char) {#1};}}

\begin{document}

\address{
Department of Mathematical Sciences\\
Norwegian University of Science and Technology\\
N-7491 Trondheim,
Norway\\
email:\,\tt{alexander.lundervold@gmail.com}
\\[4pt]
Department of Mathematics\\
University of Bergen\\
Postbox 7800, N-5020 Bergen, Norway\\
email:\,\tt{hans.munthe-kaas@math.uib.no}
}

\title{On algebraic structures of numerical integration on vector spaces and manifolds}

\author{Alexander Lundervold\thanks{Work partially supported by the
    NILS Mobility Project, UCM-EEA Abel Extraordinary Chair, 2010} and
Hans Z. Munthe-Kaas\thanks{
Work partially supported by NFR project 176891 SpadeAce}}

\maketitle

\begin{abstract} Numerical analysis of time-integration algorithms has applied advanced algebraic techniques for more than fourty years. An explicit description of the group of characters in the Butcher--Connes--Kreimer Hopf algebra first appeared in Butcher's work on composition of integration methods in 1972. In more recent years, the analysis of structure preserving algorithms, geometric integration techniques and integration algorithms on manifolds have motivated the incorporation of other algebraic structures in numerical analysis. In this paper we will survey algebraic structures that have found applications within these areas. 
This includes pre-Lie structures for the geometry of flat and torsion
free connections appearing in the analysis of numerical flows on
vector spaces. The much more recent post-Lie and D-algebras appear in the analysis of flows on manifolds with flat connections with constant torsion. Dynkin and Eulerian idempotents appear in the analysis of non-autonomous flows and in backward error analysis. Non-commutative Bell polynomials and a non-commutative Fa\`a di Bruno Hopf algebra are other examples of structures appearing naturally in the numerical analysis of integration on manifolds.
\end{abstract}



\section{Introduction}

We present some aspects of the field \emph{Geometric Numerical Integration} (GNI), whose aim it is to construct and study discrete numerical approximation methods for continous dynamical systems, exactly preserving some underlying geometric properties of the system. Many systems have conserved quantities, e.g. the energy in a conservative mechanical system or the symplectic structures of Hamiltonian systems, and numerical methods that take this into account are often superior to those constructed with the more classical goal of achieving high order.

An important tool in the study of numerical methods is \emph{Butcher series} (B-series), invented by John Butcher in the 1960s. These are formal series expansions indexed by rooted trees, and have been used extensively for order theory and the study of structure preservation. We will put particular emphasis on B-series and their generalization to methods for dynamical systems evolving on manifolds, called \emph{Lie--Butcher series} (LB-series).

It has become apparent that algebra and combinatorics can bring a lot of insight into this study. Many of the methods and concepts are inherently algebraic or combinatoric, and the tools developed in these fields can often be used to great effect. Several examples of this will be discussed throughout.

\section{Numerical integration on vector spaces}\label{sect:classgeomint}
In numerical analysis one of the main objects of study is flows of vector fields, given by initial value problems of the type\footnote{Non-autonomous differential equations can also be written on this form by adding a component to the $y$ vector}: 
\begin{equation}\label{ivp}
y'(t) = F(y(t)), \quad y(0)=y_0.
\end{equation}	
The function $y$ can be real-valued or vector-valued (giving rise to a system of coupled differential equations). The flow of the differential equation is the map $\Psi_{t,F}: \mathbb{R}^n \rightarrow \mathbb{R}^n$ defined by $y(t) = \Psi_{t,F}(y_0).$\footnote{Here we assume Lipschitz continuity of $F$ for the flow to exist and be unique.} Note that $F(y) = d/dt \Psi_{t,F}(y_0).$ In many practical settings the vector field $F$ is Hamiltonian, and their flows have several interesting geometric properties. We seek to construct good approximations to flows, where `good' can mean several different things, depending on the context. Sometimes what we want are approximations of high order, other times we need them to preserve some qualitative or geometric structure of the underlying dynamical system. Preserving geometric structure is particularly important when studying systems over long time intervals. An early illustration of this fact was made by Wisdom and Holman in \cite{wisdom1991smf}, where they computed the solar system's evolution over a billion-year time period using a symplectic method, making an energy-error of only $2 \times 10^{-11}$. 
%The term \emph{Geometric Numerical Integration} denotes a research area where emphasis is towards the construction of discrete numerical methods which exactly preserve some underlying geometric properties of the continuous dynamical system. 

As there are several excellent introductions to geometric numerical integration on vector spaces we will not go into a detailed study here, but merely describe some of the main ideas. The book \cite{hairer2006gni} is the standard reference; other introductions can be found in \cite{ mclachlan2006gif, leimkuhler2004shd, budd2003gia, mclachlan2001slo,  sanz-serna1994nhp, tsethesis, vilmart08edi}. 

The focus of this paper will be on some of the algebraic and combinatorial tools of geometric numerical integration, with particular emphasis on the tools we will use when studying flows on more general manifolds in Section \ref{sect:gengeomint}. Lately, there has been quite a lot of interest in these algebraic aspects of geometric integration, and this has resulted in both an increased understanding of the field, and also of its relations to other areas of mathematics.

\subsection{Numerical methods and structure-preservation}\label{structurepreservation}

Consider an initial value problem of the form (\ref{ivp}):
\[y'(t) = F(y(t)), \quad y(0)=y_0 \tag{1}\]
representing the flow of the (sufficiently smooth) vector field $F$. A
numerical method for (\ref{ivp}) generates approximations $y_1, y_2,
y_3, \dots$ to the solution $y(t)$ at various values of $t$. One of
the simplest methods is the (explicit) \emph{Euler method}. It
computes approximations $y_n$ to the values $y(nh)$, where $n \in \mathbb{N}$ and $h$ is the step size, using the rule:
\begin{equation}\label{classeuler}
y_{n+1} = y_n + hF(y_n).
\end{equation}
This generates a numerical flow $\Phi_h$ approximating the exact flow $\Psi$ of $F$. The accuracy of the method can be measured by its \emph{order}: we say that a one-step method $y_{n+1} = \Phi_h(y_n)$ has order $n$ if $|\Phi_h(y) - \Psi_h(y)| = O(h^{n+1})$ as $h \rightarrow 0$. Another way to put this is in terms of the curve traced out by the numerical flow: by comparing its Taylor series to the Taylor series for the curve of the exact flow term by term, we can read off the order of the method. The Taylor series for $y$ has the form
\begin{equation}
y(h) = y_0 + hF(y_0) + \frac 12 h^2 F'(y_0)F(y_0) + O(h^3),
\end{equation}
and we note that the Euler method is of order 1.

\paragraph{Runge--Kutta methods.} The Euler method is an example of a \emph{ Runge--Kutta method}, a class of methods that are extremely prevalent in applications \cite{hairer1993sod, butcher2008nmf}. A Runge--Kutta method is a one-step method computing an approximation $y_1$ to $y(h)$ with $y_0$ as input, as follows:
\begin{definition}\label{def:RK}
An $s$-stage Runge--Kutta method for solving the initial value problem (\ref{ivp}) is a one-step method given by
\begin{align}
\begin{split}
Y_i =& y_0 + h\sum_{j=1}^s a_{ij} F(Y_j), \,\,\, i = 1,\dots s \\
y_1 =& y_0 + h\sum_{i=1}^s b_i F(Y_i),
\end{split}
\end{align}
where $b_i, a_{ij} \in \mathbb{R}$, $h$ is the step size and $s \in \mathbb{N}$ denotes the number of \emph{stages}, i.e.\ the number of evaluations of $F$.
\end{definition}
A Runge--Kutta method can be presented as a \emph{Butcher tableau}, which characterizes the method completely: 
\begin{center}
\begin{tabular}{c|ccc}
$c_1$ & $a_{11}$ &$\dots$ & $a_{1s}$ \\
$\vdots$ & $\vdots$ && $\vdots$ \\
$c_s$ & $a_{s1}$ & $\dots$ &$a_{ss}$\\
\hline
& $b_1$ &$\dots$ & $b_s$
\end{tabular}
\end{center}
The coefficients $c_i = \sum_{j=1}^s a_{ij}$,  appear explicitly only  in integration methods for non-autonomous equations.
The method is called \emph{explicit} if the matrix $\{a_{i,j}\}$ is strictly lower triangular, and \emph{diagonally implicit} if it is
lower triangular. Otherwise the method is called \emph{implicit}.

\begin{example}
We note that the Euler method~(\ref{classeuler}) is the Runge--Kutta method with Butcher tableau: 
\begin{center}
\begin{tabular}{c|c}
$0$ & $0$ \\
\hline
& $1$ 
\end{tabular}
\end{center}
Another well-known example is the explicit midpoint method: 
\begin{equation}
y_{n+1} = y_n + hF\left(y_n + \frac 12 hF(y_n)\right),
\end{equation}
given by: 
\begin{center}
\begin{tabular}{c|cc}
$0$ & $0$ & $0$\\
$1/2$ & $1/2$ & $0$ \\
\hline
& $0$ &  $1$ 
\end{tabular}
\end{center}
An implicit method with good structure preserving properties is the
implicit midpoint rule 
\begin{equation}
y_{n+1} = y_n + hF\left(\frac{y_n+y_{n+1}}{2}\right),
\end{equation}
represented by the tableau
\begin{center}
\begin{tabular}{c|c}
$\frac12$ & $\frac12$ \\
\hline
& $1$ 
\end{tabular}
\end{center}
\end{example}

Given any number $m$, there exist  Runge--Kutta methods of order $m$ \cite{butcher2008nmf}. Verifying this involves expanding the methods into series involving the derivatives of $F$, and already at low orders the expressions get quite complicated. However, in Section \ref{B-series} we shall see that the Runge--Kutta methods are special cases of \emph{Butcher series methods}, and that one can find nice descriptions of the order theory and also structure preservation properties for numerical methods within this framework.

\paragraph{Differential equations and geometric structures.}

When presented with a system modeled by a differential equation one will often first try to determine its qualitative properties: are there any invariants? What kind of geometric structure does the system have? Structures of interest can be energy and volume preservation, symplectic structure, first integrals, restriction to a particular manifold (as studied in Section \ref{sect:gengeomint}), etc. Then, when choosing (or designing) a numerical method for approximating the solution of the differential equation, it might make sense for the method to share these qualitative features. In that way one has control over what kind of errors the method introduces, obtaining a method tailor-made to the problem at hand.

A rich source of problems with geometric structures are the \emph{Hamiltonian systems}. Let $H: \RR^{2n} \rightarrow \RR$ be a smooth function. A \emph{ Hamiltonian vector field} is a vector field on $\RR^{2n}$ of the form $X_H = \Omega^{-1}\nabla H$, where $\Omega$ is an antisymmetric, invertible $2n\times 2n$ matrix.\footnote{Hamiltonian vector fields can be defined on any symplectic manifold \cite{arnold1989mmo}.} %The function $H$ is called the \emph{ Hamiltonian function} of $X_H$.  
The flow of $X_H$ is given by 
\begin{equation}
\frac{d}{dt} z = \Omega\nabla_z H(z).
\end{equation}
The function $H$ represents the total energy of the system. Two important properties of the flow of a Hamiltonian vector field $X_H$ is that it is constant along the Hamiltonian function $H$  (conservation of energy) and that it preserves a symplectic form $\omega$ on $R^{2n}$. Using numerical integrators constructed to preserve these properties has been shown to lead to dramatic improvements in accuracy. For examples of this phenomenon see e.g.  \cite{hairer2006gni, hairer2005iao, leimkuhler2004shd}, and references therein.

\subsection{Trees and Butcher series}\label{B-series}

Starting with the work of John Butcher in the 1960s and 70s \cite{butcher1963cft, butcher1972aat} the study of methods for solving ordinary differential equations has been closely connected to the combinatorics of rooted trees. Many numerical methods $y_{n+1} = \Phi_h(y_n)$ (including all Runge--Kutta methods) can be expressed as certain formal series, called \emph{ Butcher series} by Hairer and Wanner in \cite{hairer1974otb}. By a clever representation of the terms the series can be indexed over the set of rooted trees. 

Consider the differential equation
\begin{equation}
y'(t) = F(y(t)).
\end{equation}
Denote the components of $F: \mathbf{R}^n \rightarrow \mathbf{R}^n$ by $f^i$ and write 
\begin{equation}\label{eq:components}
f^i_{j_1j_2\cdots j_k} = \frac{\partial^k f^i}{\partial x_{j_1} \partial x_{j_2} \cdots \partial x_{j_k}}.
\end{equation}
Summing over repeated indices, the first few derivatives of $y$ can be written as:
\begin{align}
\begin{split}
\frac{dy^i}{dt} &= f^i \\
\frac{d^2y^i}{dt^2} &= f^i_j f^j \\
\frac{d^3y^i}{dt^3} &= f^i_{jk}f^jf^k + f^i_j f^j_k f^k \\
\frac{d^4y^i}{dt^4} &= f^i_j f^j_k f^k_l f^l + f^i_j f^j_{kl} f^k f^l + 3f^i_{jk} f^j_l f^k f^l + f^i_{jkl} f^j f^k f^l.
\end{split}
\end{align}
These expressions soon get very complicated, but the structure can be made much more transparent by observing that the derivatives of $F$ can be associated in a bijective way with rooted trees, an observation already made by Cayley in 1857 \cite{cayley1857taf}. Before giving the exact correspondence between differential equations, rooted trees and Butcher series, we will take a closer look at trees.

\paragraph{Rooted trees.}\label{sect:trees}
A \emph{tree} is a connected graph with no cycles 
$$ T = \{\ab,\aabb,\aababb,\aaabbb, \aabababb, \aabaabbb, \aaababbb,\ldots \}.$$
A \emph{rooted tree} is a tree with one vertex designated as the root. In the pictorial representation of trees, the root will always be drawn as the bottom vertex, and the trees will be ordered from the root to the top. More precisely, a tree $\tau$ is a graph consisting of a set of vertices $V(\tau)$ and edges $E(\tau)\subset V(\tau) \times V(\tau)$ so that there is exactly one path connecting any two vertices.  A \emph{ path} between $v_i$ and $v_j$ is a set of edges $\{v_{s_l}, v_{t_l}\}$ so that $l=1,2,\dots, r$, $s_1=i$, $t_l = s_{l+1}$ and $t_r=j$. This gives a partial ordering of the tree in terms of paths from the root to the vertices of the tree. A vertex $v_i$ is smaller than another distinct vertex $v_j$, e.g. $v_i\prec v_j$, if the unique path from the root to $v_j$ goes via $v_i$. A vertex $v_i$ is called a \emph{ leaf} if there is no vertex $v_j$ with $v_i \prec v_j$. A \emph{child} of a vertex $v_i$ is a vertex $v_j$ with $v_i \prec v_j$ so that there is no vertex $v_k$ with $v_i\prec v_k \prec v_j$. The \emph{ order} $|\tau|$ of a tree $\tau$ is the number of vertices of the tree. We define a symmetry group on a tree $\tau$ as all automorphisms on the vertices. The order of this group, $\sigma(\tau)$, is called the \emph{symmetry} of the tree $\tau$. For example, 
\[\sigma\left(\aabb\right) = 1, \quad \sigma\left(\aababb \right) = 2, \quad \sigma \left( \aabababb\right) = 6.\]
A recursive definition of $\sigma$ can be found in \cite{hairer2006gni}.

A \emph{forest} of rooted trees is a graph whose connected components are rooted trees, e.g. $\omega = \tau_1 \dots \tau_n.$ We include the \emph{empty tree} $\one$, i.e. the graph with no vertices, in the set of forests. The set of forests can be put in bijection to the set of trees via the operator $B^+:F \rightarrow T$, defined on a forest $\omega = \tau_1 \dots \tau_n$ by connecting the trees to a new root by addition of edges. For example, 
\[B^+(\aababb\, \ab) = \aaababbabb.\]
This operator can be used to generate all trees recursively from the tree $\ab$ by the following procedure:

\begin{itemize}
\item[(i)] The graph $\ab$ belongs to $\T$
\item[(ii)] If $\tau_1,\dots, \tau_n \in \T$ then $\tau = B^+(\tau_1\dots \tau_n)$ is in $\T$. 
\end{itemize}

\noindent The \emph{tree factorial} $\tau!$ is given recursively by:
\begin{itemize}
\item[(i)] $\ab! = 1$ 
\item[(ii)] $B^+(\tau_1 \dots \tau_n)! = |B^+(\tau_1 \dots \tau_n)|\tau_1! \dots \tau_n!$.
\end{itemize}
%Another important number associated to a tree is the \emph{Connes--Moscovici coefficient} \cite{kreimer1999cii, brouder2000rkm}:
%\begin{equation}
%CM(\tau) = \frac{|\tau|!}{\tau ! \,\sigma(\tau)}.
%\end{equation}
%We shall see that these show up in the Butcher series for the exact
%flow of the initial value problem (\ref{ivp}). 
An important operation on trees is the Butcher product, defined in terms of \emph{grafting}. 
\begin{definition}\label{bpr}
The \emph{Butcher product} $\tau \bpr \omega$ of a tree $\tau = \Bplus(\tau_1\dots \tau_n)$ and a forest $\omega = \omega_1 \cdots \omega_m$ is given by grafting $\omega$ onto the root of $\tau$: 
\begin{equation}
\tau \bpr \omega = \Bplus(\tau_1 \dots \tau_n\, \omega_1\dots \omega_m) 
\end{equation}
\end{definition}

\paragraph{Butcher series.}
The calculations of the derivatives of $y'(t) = F(y(t))$ performed at the beginning of the section can be written in terms of the elementary differentials of $F$.
\begin{definition}\label{def:elementdiff}
Let $F:\mathbb{R}^n \rightarrow \mathbb{R}^n$ be a vector field. The \emph{elementary differential} $\mathcal{F}$ of $F$ is
\begin{align}
\begin{split}
\mathcal{F}(\ab)(y) &= F(y) \\
\mathcal{F}(\tau)(y) &= F^{(m)}(y)(\mathcal{F}(\tau_1)(y), \dots, \mathcal{F}(\tau_m)(y)),
\end{split}
\end{align}
where $F^{(m)}$ is the $m$-th derivative of the vector field $F$ and $\tau = B^+(\tau_1, \dots, \tau_m)$ is a rooted tree. 
\end{definition}

\noindent We will discuss another way to write elementary differentials in Section \ref{sect:prelie}. With the notation from Equation (\ref{eq:components}), the first few elementary differentials are shown in Table (1.1). The vector field $F$ corresponds to the leaves of the tree, the first derivative $F'$ corresponds to a vertex with an edge with one child,  the second derivative $F''$ corresponds to a vertex with two children, etc.

\begin{table}[h!]~\label{table:elementdiff}
  \begin{equation*}
    \begin{array}{c|c}
      \tau & \F(\tau)(y)^i \\[1mm]
      \hline\\[-2mm]
      \ab& f^i \\[1mm]
      \aabb & f^i_jf^j \\[1mm]
      \aababb & f^i_{jk}f^jf^k \\[1mm]
      \aaabbb & f^i_jf^j_kf^k \\[1mm]
      \aabababb & f^i_{jkl}f^jf^kf^l \\[1mm]
      \aabaabbb & f^i_{jk}f^jf^k_lf^l 
    \end{array}
  \end{equation*}
\caption{Elementary differentials associated to a vector field $F$ with components $f^i$.}
\end{table}

Butcher series are (formal) Taylor expansions of elementary differentials indexed over trees:
\begin{definition}
A \emph{Butcher series} (B-series) is a (formal) series expansion in a parameter $h$:  
\begin{align}
\begin{split}
\Bs_{h,F} (\alpha) &= \alpha(\one)\F(\one) + \sum_{\tau \in \T}h^{|\tau|}\frac{\alpha(\tau)}{\sigma(\tau)} \F(\tau)\\ 
&= \sum_{\tau \in \tilde{\T}}h^{|\tau|}\frac{\alpha(\tau)}{\sigma(\tau)} \F(\tau),
\end{split}
\end{align}
where $\tilde{\T} = \T \cup \{\one\}$, $F$ is a vector field, $\alpha$ is a function $\alpha: \tilde{\T} \rightarrow \RR$, $\sigma(\tau)$ is the symmetry of $\tau$, $h$ is a real number (representing the step size), and $\F$ is the elementary differential of $F$, extended to the empty tree $\one$ by $\F(\one)(y) = y$.
\end{definition}

We shall see that these series can be used to represent numerical methods $y_{n+1} = \Phi_{h}(y_n)$ approximating the flow of a vector field $F$, in the sense that the Taylor series for $\Phi_h$ can be expanded into a B-series: $\Phi_{h} = \Bs_{h,F}(\alpha)$.\footnote{A numerical method for solving a differential equation is called a \emph{B-series method} if it can be written as a B-series.}

By computing the Taylor expansion of the solution to the initial value problem (\ref{ivp}) one obtains the following result:
\begin{proposition}[{\cite{hairer2006gni}}]
The Taylor series for the solution of the differential equation (\ref{ivp}) can be written as a B-series:
\begin{equation}
B_{h,F}(\gamma) = \sum_{\tau \in \tilde{\T}} h^{|\tau |} \frac{\gamma(\tau)}{\sigma(\tau)} \mathcal{F}(\tau),
\end{equation}
where $\gamma(\tau) = 1/\tau!$. That is, $y(t + h) = \Bs_{h,F}(\gamma)(y(t))$. 
\end{proposition}

Runge--Kutta methods can also be written as B-series expansions, with coefficients given by the \emph{elementary weights} of the method \cite{butcher1963cft}.

\begin{definition}[Elementary weights]
Let $b_i$ and $a_{ij}$ be coefficients of a RK-method as in Definition \ref{def:RK}, where $i \in \mathbb{N}$. The \emph{elementary weight function} $\Phi$ is defined on trees as follows:
\begin{align}
\begin{split}
&\Phi_i(\ab) = c_i\\
&\Phi(\ab) = \sum_{j=1}^s b_j \\
&\Phi_i(B^+(\tau_1,\dots, \tau_k)) = \sum_{j=1}^s a_{ij} \Phi_j(\tau_1) \Phi_j(\tau_2) \dots \Phi_j(\tau_k)\\
&\Phi(B^+(\tau_1,\dots, \tau_k)) = \sum_{j=1}^s b_j \Phi_j(\tau_1) \Phi_j(\tau_2) \dots \Phi_j(\tau_k)  
\end{split}
\end{align}
Here $i=1, \dots, s$.
\end{definition}
For example,
\[\Phi(\aabb) = \sum_{j=1}^s b_jc_j, \quad \Phi(\aababb) = \sum_{j=1}^s b_jc_j^2, \quad \Phi(\aaababbb) =  \sum_{j,k = 1}^s b_ja_{jk}c_k^2\]
\begin{theorem}[{\cite{butcher1963cft}}] The B-series for a RK-method given by the elementary weights $\Phi(\tau)$ is
\begin{equation}
\Bs_{h,F}(\Phi) = \sum_{\tau \in \tilde{\T}} h^{|\tau|}\frac{\Phi(\tau)}{\sigma(\tau)}\F(\tau)
\end{equation}
\end{theorem}

\paragraph{Order theory for B-series methods.}
Once we have the B-series of the exact solution and the B-series of a numerical method, it is straightforward to compare the coefficients and read off the order of the method. For Runge--Kutta methods, we obtain the following result:
\begin{proposition}[{\cite{butcher1963cft}}]
A Runge--Kutta method given by a B-series with coefficients $\Phi(\tau)$ has order $n$ if and only if 
\begin{equation}
\Phi(\tau) = \gamma(\tau), \hspace{0.5cm}\text{  for all } \tau \in T
\text{  such that } |\tau| < n.
\end{equation}
\end{proposition}

\paragraph{B-series methods and structure preservation.} The class of
B-series methods includes all Taylor series methods and Runge--Kutta
methods. It does not, however, include all numerical methods, an
example being the class of \emph{splitting methods}.

It is important to point out that focusing only on B-series methods has its drawbacks. Besides the fact that the class does not contain all methods, it is also known that there are certain geometric structures that cannot be preserved by B-series methods. For example, no B-series method can preserve the volume for all systems \cite{iserles2007bsm}. However, we will be content with this loss of generality and focus exclusively on methods based on B-series in this section, and on their generalization -- Lie--Butcher series -- in the next.

A case which is particularly well-studied is Hamiltonian vector fields. The following two theorems serve as prime examples:

\begin{theorem}[\cite{hairer1994bao}] Let $G = \Bs_{h,F}(\alpha)$ be a vector field with $\alpha(\one) = 0,$ $\alpha(\ab) \neq 0$. Then $G$ is Hamiltonian for all Hamiltonian vector fields $F(y) = \Omega^{-1} \nabla H(y)$ if and only if
\begin{equation}
\alpha(\tau_1 \bpr \tau_2) + \alpha(\tau_2 \bpr \tau_1) = 0
\end{equation}
for all $\tau_1, \tau_2 \in \T$. Here $\bpr$ denotes the Butcher product of Definition \ref{bpr}.
\end{theorem}

\begin{theorem}[\cite{calvo1994cbs}] Consider a numerical method given by a B-series $\Bs_{h,F}(\alpha)$. The method is symplectic if and only if 
\begin{equation}
\alpha(\tau_1 \bpr \tau_2) + \alpha(\tau_2 \bpr \tau_1) = \alpha(\tau_1)\alpha(\tau_2)
\end{equation}
for all $\tau_1, \tau_2 \in \T$, where $\alpha(\one)=0$.
\end{theorem}

The paper \cite{celledoni2010epi} gives an overview of what is known about structure preservation for B-series, including characterizations of the various subsets of trees corresponding to energy-preserving, Hamiltonian and symplectic B-series.

\subsection{Hopf algebras and the composition of Butcher series}\label{GI:composition}

Consider two numerical methods given by $\Phi^1$ and $\Phi^2$. Using the method $\Phi^1$ to advance a point $y_0$ to a point $y_1$, and then applying the method $\Phi^2$ using $y_1$ as initial point, results in a point $y_2$:
\[y_1 = \Phi^1(y_0), \hspace{1cm} y_2 = \Phi^2(y_1).\]
This is the idea behind \emph{composition} of numerical methods. In the case where both methods are given by B-series, $\Phi^1 = \Bs^1_{h,F}(\alpha)$, $\Phi^2 = \Bs^2_{h,F}(\beta)$, then the composition method $\Phi^2 \circ \Phi^1$ is again a B-series: $\Phi^2\circ \Phi^1 = \Bs_{h,F}(\gamma)$. Concretely, if $y_1 = \Bs^1_{h,F}(\alpha)(y_0)$ and $y_2 = \Bs^2_{h,F}(\beta)(y_1)$, then $y_2 = \Bs_{h,F}(\gamma)(y_0)$. This is the Hairer--Wanner theorem from \cite{hairer1974otb}. The coefficient function $\gamma$ of this B-series was first studied by John Butcher in \cite{butcher1972aat}, where he found that composition of B-series is a group operation (giving rise to the \emph{Butcher group}) on the coefficient functions, and gave expressions for the product, identity and inverse in this group.

In \cite{kreimer1998oth, connes1998har} Connes and Kreimer introduced a Hopf algebra of rooted trees connected to the renormalization procedure in  quantum field theory. Later \cite{brouder2000rkm} it was pointed out that a variant of this Hopf algebra is closely related to the Butcher group. More precisely, the Butcher group is the group of characters in a Hopf algebra $\BCK$ defined by Connes and Kreimer.

We will describe the Butcher group indirectly by describing the Hopf algebra $\BCK$. But first we will present some basic definitions from the theory of Hopf algebras. For a comprehensive introduction, see \cite{sweedler1969ha, abe80ha}. Other excellent references include \cite{cartier2006apo, manchon2008hai}.

\paragraph{Hopf algebras.} 

Let $\field$ be a field of characteristic zero. An \emph{algebra} $A$ over $\field$ is a $\field$-vector space equipped with a multiplication map $\mu: A \otimes A \rightarrow A$ and a unit $u: \field \rightarrow A$ so that 
\begin{itemize}
\item \begin{tabular}{lr} $\mu \circ (id \otimes \mu) = \mu \circ (\mu \otimes id): A\otimes A \otimes A \rightarrow A$ & (associativity)\end{tabular}
\item \begin{tabular}{lr} $\mu \circ (u \otimes id) = \mu \circ (id \otimes u): k \otimes A \cong A \rightarrow A$ &(unitality)\end{tabular}
\end{itemize}
A \emph{coalgebra} $C$ over $\field$ is the dual notion. It consists of a comultiplication map $\Delta: C \rightarrow C \otimes C$ and a counit $\epsilon: C \rightarrow \field$ so that
\begin{itemize}
\item \begin{tabular}{lr} $(\Delta \otimes id) \circ \Delta = (id \otimes \Delta) \circ \Delta: C \rightarrow C \otimes C \otimes C$ & (coassociativity)\end{tabular}
\item \begin{tabular}{lr} $(\epsilon \otimes id) \circ \Delta = (id \otimes \epsilon) \circ \Delta: C \rightarrow C \otimes \field \cong C$ &(counitality)\end{tabular}
\end{itemize}
A \emph{Hopf algebra} is at once an algebra and a coalgebra, and it comes equipped with an antipode $S: H \rightarrow H$. These structures have to satisfy certain compatibility conditions, written as the following diagrams, where $\tau$ denotes the flip operation $\tau(h_1, h_2) = (h_2, h_1)$:

\begin{minipage}[b]{0.5\linewidth}
\begin{center}
\begin{diagram}[labelstyle=\scriptstyle]
H^{\otimes 4} &&\rTo^{I \otimes \tau \otimes I}& &  H^{\otimes 4}\\
\uTo^{\Delta \otimes \Delta}    &&&&  \dTo_{\mu \otimes \mu}\\
H \otimes H&\rTo_{\mu}& H &\rTo_{\Delta}& H \otimes H
\end{diagram}
\end{center}
\end{minipage}
\begin{minipage}[b]{0.5\linewidth}
\begin{center}
\begin{diagram}[labelstyle=\scriptstyle]
H \otimes H &\rTo^{\epsilon \otimes \epsilon}& k \otimes k\\ 
\dTo^{\mu} && \dTo_{\cong} \\
H &\rTo_{\epsilon}& k
\end{diagram}
\end{center}
\end{minipage}
\vspace{0.5cm}
\begin{diagram}
& & H\otimes H & & \rTo^{\scriptstyle S\otimes 1} & & H\otimes H \\
& \ruTo^{\scriptstyle\Delta} & & & & & & \rdTo>{\scriptstyle\mu} \\
H & & \rTo^{\scriptstyle \varepsilon} & & k & & \rTo^{\scriptstyle u} & & H \\
& \rdTo<{\scriptstyle\Delta} & & & & & & \ruTo>{\scriptstyle \mu} \\
& & H\otimes H & & \rTo_{\scriptstyle 1\otimes S} & & H\otimes H \\
\end{diagram}
The first two diagrams ensure that the coproduct and the counit are both algebra homomorphisms. The last diagram is best interpreted in terms of the characters in a Hopf algebra. Let $A$ be a commutative $k$-algebra, and let $\Lin(H,A)$ denote the set of linear maps from $H$ to $A$. An element $\alpha \in \Lin(H,A)$ is called a \emph{character} if $\alpha(x\cdot y) = \alpha(x) \cdot \alpha(y)$ for all $x,y \in H$, where the product on the left-hand side is in $H$, and on the right-hand side in $A$. The set of characters in $\Lin(H,A)$ form a group under the \emph{convolution product}: 
\begin{equation}\label{eq:convolution}
\phi \star \psi = \mu \circ (\phi \otimes \psi) \circ \Delta.
\end{equation}
The unit is the composition of the unit and the counit in $H$,
e.g. $\eta := u \circ \epsilon.$ The bottom diagram above corresponds
to the antipode being the inverse of the identity under this product,
and we have $\alpha^{\star -1}=\alpha \circ S$.

Later we will also need the concept of \emph{infinitesimal characters} (also called \emph{derivations}),
which are maps $\alpha$ in $\Lin(H,A)$ satisfying 
\begin{equation}
\alpha(x\cdot y) = \eta(x)\cdot\alpha(y) + \alpha(x)\cdot\eta(y).
\end{equation}

A Hopf algebra $H$ is \emph{graded} if it is graded as an algebra, i.e. $H = \bigoplus_{n\geq 0} H_n$ with $\mu(x_r, x_s) \in H_{r+s}$ for $x_r \in H_r$ and $x_s \in H_s$, and its coproduct satisfies 
\[\Delta(H_n) \subset \bigoplus_{r+s = n} H_r \otimes H_s.\]

\paragraph{The Butcher--Connes--Kreimer Hopf algebra.} The composition of B-series is governed by a certain Hopf algebra $\BCK$ based on the set $T$ of rooted trees, called the \emph{Butcher-Connes-Kreimer Hopf algebra}. In Section \ref{sect:gengeomint} we will see that a generalization of this Hopf algebra governs the composition of Lie-Butcher series (Section 2.\ref{sect:compLB}). 

To describe the BCK Hopf algebra we need to define its structure as a vector space, an algebra, a coalgebra, and define the antipode. As a $\RR$-vector space $\BCK$ is generated by the set $T$ of rooted trees, and graded by the order (i.e. number of vertices) of the trees. The algebra structure is that of the symmetric algebra $S(\RR\{T\})$. The product is written as (commutative) concatenation of trees (i.e. disjoint union), giving rise to forests of trees. The unit is the empty tree $\one$.
\[\aabb \, \aababb = \aababb \, \aabb, \qquad \aabb \,\,\, \one = \one \,\,\,\aabb = \aabb\]
The coproduct of $\BCK$ is the map $\DeltaBCK: \BCK \rightarrow \BCK \otimes \BCK$ determined recursively by:
\begin{equation}
\DeltaBCK \circ B^+(\omega) = B^+(\omega) \otimes \one + (Id \otimes B^+) \circ \DeltaBCK(\omega),
\end{equation}
where $\omega$ is a forest\footnote{Recall that $\DeltaBCK$ is an algebra morphism and is therefore defined on forests as well as trees, since $\DeltaBCK(\tau_1\tau_2) = \DeltaBCK(\tau_1)\DeltaBCK(\tau_2)$.}. The counit is the map $\epsilon: \BCK \rightarrow  \RR$ given by $\epsilon(\one)= 1$ and $\epsilon(\tau) = 0$ if $\tau \neq \one$. The coproduct can also be written in a non-recursive manner using cuttings of trees.

\paragraph{Cutting trees.} 
An \emph{admissible cut} of a tree $\tau$ is a set $c \subset E(\tau)$ of edges of $\tau$ such that $c$ contains at most one edge from any path from the root to a leaf. The case $c = \emptyset$ is called the empty cut. Let $\omega$ denote the forest with vertices $V(\tau)$ and edges $E(\tau)\setminus c$. We write $R^c(\tau)$ for the component of $\omega$ containing the root of $\tau$, and $P^c(\tau)$ for the forest consisting of the remaining components. The cut resulting in $P^c(\tau) = \tau$ and $R^c(\tau) = \one$ is also admissible, and called the \emph{full cut} (f.c.).

\begin{theorem}[{\cite{connes1998har}}]\label{thm:BCKcoprod}
The coproduct in $\BCK$ can be written as 
\begin{equation}
\DeltaBCK(\tau) = \sum_{c\in Adm(\tau)} P^c(\tau) \otimes R^c(\tau)
\end{equation}
\end{theorem}
\noindent Examples of the coproduct can be found in Table \ref{BCKcoprod}.
\begin{table}[h! ]
  \centering
  \begin{equation*}
    \begin{array}{c@{\,\,}|@{\quad}l}
      \hline \\[-2mm]
      \tau & \DeltaBCK(\tau)  \\[1mm]
      \hline \\[-2mm]
      \one & \one\tpr\one \\[1mm]
      \ab & \ab\tpr\one+\one\tpr\ab \\[2mm]
      \aabb & \aabb\tpr\one+\ab\tpr\ab+\one\tpr\aabb  \\[2mm]
      \aaabbb &\aaabbb\tpr\one+\ab\tpr\aabb+\aabb\tpr\ab+\one\tpr\aaabbb \\[2.5mm]
      \aababb & \aababb\tpr\one+\ab\ab\tpr\ab+2\,\ab\tpr\aabb+\one\tpr\aababb \\[2.5mm]
      \aaaabbbb & \aaaabbbb\tpr\one+\aaabbb\tpr\ab+\aabb\tpr\aabb+
      \ab\tpr\aaabbb+\one\tpr\aaaabbbb \\[2.5mm]
      \aaababbb & \aaababbb\tpr\one+\aababb\tpr\ab+\ab\ab\tpr\aabb+
      2\,\ab\tpr\aaabbb+\one\tpr\aaababbb \\[2.5mm]
      \aabaabbb & \aabaabbb \otimes \one + \ab\aabb\otimes \ab + \ab\ab \otimes \aabb + \aabb\otimes \aabb + \ab \tpr \aababb + \ab \otimes \aaabbb +  \one \otimes \aabaabbb
    \end{array}
  \end{equation*}
  \caption{Examples of the coproduct $\DeltaBCK$ in the Hopf algebra $\BCK$}\label{BCKcoprod}
\end{table}
The antipode can be defined recursively as $S(\one) = \one$ and: 
\begin{equation}
S(\tau) = -\tau - \sum _{c \in Adm(\tau)\setminus \{\emptyset \cup f.c.\}}S(P^c(\tau))R^c(\tau)
\end{equation}

The Hairer--Wanner theorem gives the exact correspondence between $\BCK$ and composition of B-series:

\begin{theorem}[{\cite{hairer1974otb}}]\label{hairerwannertheorem}
Let $\Bs^1_{h,F}(\alpha)$ and $\Bs^2_{h,F}(\beta)$ be two B-series, with coefficients $\alpha, \beta: T \rightarrow \RR$. The composition $\Bs^2_{h,F}(\beta) \circ \Bs^1_{h,F}(\alpha)$ is again a B-series, and we have 
\begin{equation}
\Bs^2_{h,F}(\beta) \circ \Bs^1_{h,F}(\alpha) = \Bs_{h,F}(\alpha \star \beta),
\end{equation}
where $\star$ denotes convolution in the Hopf algebra $\BCK$.
\end{theorem}

\subsection{Substitution and backward error analysis for Butcher series}\label{substB}

Consider a numerical method $\Phi_h$ used to solve a differential equation of the form
\begin{equation}\label{diffeq}
y' = F(y).
\end{equation}
The basic idea of \emph{backward error analysis} of the method $\Phi_h$ is to interpret it as giving the exact solution of a modified equation:
\begin{equation}\label{modifieddiffeq}
\tilde{y}' = \tilde{F_h}(\tilde{y}).
\end{equation}
If we can find such an equation, we can use it to study the properties
of the numerical method. In other words, the numerical method $\Phi_h$
will be represented by a modified vector field $\tilde{F}$, which then can be used to study the method. The idea is based on work by Wilkinson in the context of algorithms for solving equations given by matrices \cite{wilkinson1960eao}, and has been explored in several papers \cite{warming1974tme, hairer1994bao, calvo1994mef, hairer2006gni, chartier2007nib}. Recurrence formulas for the modified equation were first obtained in \cite{hairer1994bao, calvo1994mef}.

A related notion is the \emph{modifying integrators} of \cite{chartier2007nib}. The idea is to look for a vector field $\tilde{F}_h$ so that the numerical method $\Phi_h$ applied to the flow equation of $\tilde{F}_h$ (Equation \ref{modifieddiffeq}) is the exact solution of Equation \ref{diffeq}.

It turns out that the case where $\Phi_h$ is a B-series method is particularly nice \cite{chartier2005asl, chartier2007nib, calaque2009tih}. The vector fields $\tilde{F_h}$ can then be written as B-series whose coefficients are derived from the coefficients of $\Phi_h$, and these coefficients can be expressed by the \emph{substitution law} for B-series methods (Corollaries \ref{corBEA} and \ref{corMI}).

\paragraph{The substitution law.} Let $\Bs_{h,F}(\alpha)$ and $\Bs_{h,G}(\beta)$ be two B-series, where $\alpha(\one) = 0$. Then $\Bs_{h,F}(\alpha)$ is a vector field, and we can consider the B-series obtained by using this as the vector field $G$ in the B-series $\Bs_{h,G}(\beta)$. This is called \emph{substitution of B-series}. The result is given in terms of a bialgebra $\CEFM$ by the following theorem:

\begin{theorem}[{\cite{chartier2005asl, chartier2007nib, calaque2009tih}}]\label{thm:substB}
Let $F$ be a vector field, $\alpha, \beta$ linear maps $\alpha,\beta: \T \rightarrow \RR$ where $\beta$ is an infinitesimal character of $\BCK$, and $\alpha(\one) = 0$. Then the vector field $(1/h)\Bs_{h,F}(\alpha)$ inserted into the B-series $\Bs_{h,\cdot}(\beta)$ is again a B-series, given by
\begin{equation}
\Bs_{h, (1/h) \Bs_{h,F}(\alpha)}(\beta) = \Bs_{h,F} (\alpha \star \beta),
\end{equation}
where $\star$ denotes convolution of characters in the bialgebra $H_{CEFM}$.
\end{theorem}
The bialgebra $\CEFM$ is the symmetric algebra over rooted trees $S(\T)$, with $\ab$ as unit, equipped with a coproduct given by contracting subforests in trees:
\begin{equation}\label{CEFMcoprod}
\Delta(\tau) = \sum_{\omega \subseteq \tau} \omega \tpr \tau/\omega.
\end{equation}
If $\tau$ is a tree then the notation $\omega \subset \tau$ means that $\omega$ is a spanning subforest of $\tau$, i.e. that $\omega$ is a collection of subtrees of $\tau$ so that each vertex of $\tau$ belongs to exactly one tree in $\omega$. Then $\tau/\omega$ denotes the tree obtained by contracting each subtree (with at least two vertices) of $\tau$ contained in $\omega$ onto a vertex. Some examples of the coproduct can be found in Table \ref{table:CEFM}. The bialgebra is graded by the number of edges. 

There is a Hopf algebra related to $\CEFM$, obtained by considering the symmetric algebra over the set of rooted trees $\T'$ \emph{with at least one edge} (e.g. $\ab$ is not included), and then adding $\ab$ back as the \emph{unit} for the product. The coproduct is defined as in Equation (\ref{CEFMcoprod}). The resulting bialgebra is \emph{connected}, which makes it a Hopf algebra \cite{manchon2008hai}.

For details on these constructions, consult \cite{calaque2009tih}.
\begin{table}[!ht]
  \centering
  \begin{equation*}
    \begin{array}{c@{\,\,}|@{\quad}l}
      \hline \\[-2mm]
      \tau & \Delta_{CEFM}(\tau)  \\[1mm]
      \hline \\[-2mm]
      \ab & \ab\tpr\ab \\[2mm]
      \aabb & \aabb\tpr\ab+\ab\ab\tpr\aabb  \\[2mm]
      \aaabbb &\aaabbb\tpr\ab+\ab\ab\ab\tpr\aaabbb+2\,\aabb\ab\tpr\aabb \\[2.5mm]
      \aababb & \aababb \tpr \ab + \ab\ab\ab\tpr \aababb + 2\,\aabb\ab\tpr \aabb \\[2.5mm]
      \aaaabbbb & \aaaabbbb\tpr\ab+\ab\ab\ab\ab\tpr\aaaabbbb+2\,\aaabbb\ab\tpr\aabb+
      3\,\aabb\ab\ab\tpr\aaabbb+\aabb\aabb\tpr\aabb \\[2.5mm]
      \aaababbb &  \aaababbb\tpr \ab + \ab\ab\ab\ab\tpr \aaababbb + 2\, \aabb\ab\ab\tpr \aaabbb + \aabb\ab\ab\tpr\aababb + 2\, \aaabbb\ab\tpr \aabb + \aababb\ab\tpr \aabb\\[2.5mm]
      \aabaabbb & \aabaabbb \tpr \ab + \ab\ab\ab\ab\tpr \aabaabbb + \aabb\ab\ab\tpr \aaabbb + 2\, \aabb\ab\ab\tpr \aababb + \aabb\,\aabb \tpr \aabb + \aababb\ab\tpr \aabb + \aaabbb\ab\tpr \aabb
    \end{array}
  \end{equation*}
  \caption{Examples of the coproduct $\Delta_{CEFM}$ in the substitution bialgebra}\label{table:CEFM}
\end{table}

\paragraph{Backward error analysis and modifying integrators.} Once Theorem \ref{thm:substB} is established one can obtain expressions for backward error analysis and modifying integrators. 
\begin{corollary}[{\cite{chartier2005asl, chartier2007nib}}) (Backward error analysis]\label{corBEA}
Let $\Bs_{G}(\gamma)$ denote the B-series for the exact flow of the vector field $G$, and let $\Bs_{F}(\alpha)$ be a B-series giving a numerical flow for $F$. The modified vector field $\tilde{F}$ given by $\Bs_{\tilde{F}}(\gamma) = \Bs_{F}(\alpha)$ is a B-series $\Bs_{F}(\beta)$ with coefficients given by 
\begin{equation}
\beta \star \gamma = \alpha 
\end{equation}
\end{corollary}
\begin{corollary}[{\cite{chartier2005asl, chartier2007nib}}) (Modifying integrators]\label{corMI}
Let $\Bs_{G}(\gamma)$ denote the B-series for the exact flow of the vector field $G$, and let $\Bs_{F}(\alpha)$ be a B-series giving a numerical flow for $F$. The modified vector field $\tilde{F}$ so that $\Bs_{\tilde{F}}(\alpha) = \Bs_{F}(\gamma)$ is a B-series $\Bs_{F}(\beta)$ whose coefficients are given by
\begin{equation}
\beta \star \alpha = \gamma 
\end{equation}
\end{corollary}

\subsection{Pre-Lie Butcher series}\label{sect:prelie}

The space of vector fields on $\RR^n$ has the structure of a \emph{ pre-Lie
  algebra}, and in this section we will see that B-series can be formulated purely in terms of this pre-Lie
structure. This allows us to lift the concept of B-series to the free pre-Lie
algebra, giving rise to \emph{ pre-Lie B-series}
\cite{ebrahimi-fard2010plb}. Viewing B-series as objects in the free
pre-Lie algebra gives a clearer focus on the core algebraic structures
at play, and it also enables the application of tools and results from
other fields where pre-Lie algebras appear. Examples of this
phenomenon can be found in \cite{ebrahimi-fard2009ama} (see Remark
\ref{remark:magnus}), \cite{calaque2009tih} and \cite{chapoton2002rta}. We give the basic constructions here because formulating Butcher series in terms of pre-Lie algebras will find an analogue in Section \ref{sect:gengeomint}, where Lie--Butcher series will be constructed from the so-called post-Lie algebras.

\paragraph{Pre-Lie algebras.} The concept of pre-Lie algebras is a relaxation of associative algebras that still preserve their \emph{Lie admissible} property. In other words, for an associative algebra $(A,*)$ antisymmetrization of the product $*$ gives a Lie bracket, making it a Lie algebra: $[a,b] = a*b - b*a$, and this property also holds for pre-Lie algebras. Note, however, that not all pre-Lie algebras are associative. They were first introduced and studied by Vinberg \cite{vinberg1963chc}, Gerstenhaber \cite{gerstenhaber1963tcs}, and Agrachev and Gamkrelidze \cite{agrachev1981caa}, under various names. A nice introduction to pre-Lie algebras can be found in \cite{manchonass}.

\begin{definition}\label{preLie}
A (left) \emph{ pre-Lie algebra}\footnote{Also called a \textit{Vinberg}, \textit{left-symmetric} or \textit{chronological} algebra} $(A, \tr)$ is a $k$-vector space $A$ equipped with an operation $\tr: A \otimes A \rightarrow A$ subject to the following relation:
\begin{align}
a_{\tr}(x,y,z) - a_{\tr}(y,x,z) = 0,
\end{align}
where $a_{\tr}(x,y,z)$ is the associator $a_{\tr}(x,y,z) = x \tr (y \tr z) - (x\tr y) \tr z$.
\end{definition}
\begin{example}[The pre-Lie algebra of vector fields] The space of vector fields $\mathcal{X}(M)$ on a differentiable manifold $M$ equipped with a flat, torsion-free connection $\nabla$ can be given the structure of a pre-Lie algebra by defining $\tr$ as $F \tr G = \nabla_F G$. In the case $M = \mathbb{R}^n$ with the standard flat and torsion-free connection we have that for $F = \sum_{i=1}^n F_i \partial_i$ and $G=\sum_{j=1}^n G_j \partial_j$, 
\begin{equation}
F \tr G = \sum_{i=1}^n \left(\sum_{j=1}^n F_j(\partial_jG_i)\right)\partial_i.
\end{equation}
In Section \ref{sect:gengeomint} we will see that allowing for torsion leads to the concept of \emph{ post-Lie algebras}. See also \cite{munthe-kaas2012opl}.
\end{example}

\paragraph{The free pre-Lie algebra.}  The free pre-Lie algebra has been studied in several papers, most notably by Chapoton and Livernet in \cite{chapoton2001pla}, Dan Segal in \cite{segal1994fls}, Agrachev and Gramkrelidze in \cite{agrachev1981caa}, Dzhumadil'daev and L\"{o}fwall in \cite{dzhumadildaev2002tfr}. These papers give different bases for the free pre-Lie algebra, and one can choose to work in the basis most beneficial for the problem at hand.  
A basis for the free pre-Lie algebra $PL(V)$ over a vector space $V$ was described by Chapoton and Livernet in terms of nonplanar rooted trees \cite{chapoton2001pla, chapoton2009art}:
\[\left\{\ab,\aabb,\aababb,\aaabbb, \aabababb, \aabaabbb, \aaababbb,\ldots \right\}\]
decorated by elements of $V$. The pre-Lie product $\tau_1 \car \tau_2$ of two rooted trees is given by grafting: $\tau_1 \car \tau_2$ is the sum of all the trees resulting from the addition of an edge from the root of $\tau_1$ to one of the vertices of $\tau_2$:
\begin{equation}
\tau_1 \car \tau_2 := \sum_{v \in V(\tau_2)}\tau_1 \circ_v \tau_2
\end{equation}
Here $\tau_1 \circ_{v} \tau_2$ denotes grafting at the vertex $v$ of $\tau_2$.

\[\ab \car \ab = \aabb, \qquad \ab \car \aabb = \aababb + \aaabbb, \qquad \aabb \car \aabb = \aaabbabb + \aaaabbbb\]

\begin{theorem}[{\cite{chapoton2001pla}}]\label{fpla}
$PL(V)$ is the free pre-Lie algebra on the vector space $V$: for any pre-Lie algebra $P$ equipped with a morphism $V \rightarrow
P$, there is a unique pre-Lie morphism $PL(V) \rightarrow P$ making the
following diagram commute:
\begin{diagram}[labelstyle=\scriptstyle]
V &\rTo& PL(V)\\
&\rdTo& \dTo_{\exists !}\\
&& P
\end{diagram}
\end{theorem}
We write $PL$ for the free pre-Lie algebra on a space with only one element. 

The free pre-Lie algebra is related to the Hopf algebra $\BCK$ defined in Section \ref{GI:composition}:
\begin{theorem}[{\cite{chapoton2001pla}}]
The universal enveloping algebra $U(PL)$ of the free pre-Lie algebra on the one-vertex tree, viewed as a Lie algebra, is isomorphic to the dual of the Butcher--Connes--Kreimer Hopf algebra $\BCK$.
\end{theorem}
In fact, the dual of the Butcher--Connes--Kreimer Hopf algebra is isomorphic to
the \emph{Grossman-Larson Hopf algebra} defined in \cite{grossman89has}. The isomorphism was proven in \cite{hoffman2003cor}.

\paragraph{Pre-Lie Butcher series.} Now we can formulate the pre-Lie Butcher series
\begin{definition} A \emph{ pre-Lie Butcher series} is a formal series in $\RR\langle \PL\rangle$: 
\begin{equation}\label{eq:plB}
X(\alpha) = \sum_{t \in \PL} h^{|t|}\alpha(t)t.
\end{equation}
\end{definition}
\noindent The classical B-series are recovered by applying the unique pre-Lie morphism associated to a vector field $F$: 
\begin{equation}
\F: \PL \rightarrow \mathcal{X}(\RR^n) \quad \text{such that} \quad
\F(\ab) = F.
\end{equation}
This is the elementary differential function of $F$ as defined in \ref{def:elementdiff}. It is given recursively by $\F(\ab) = F$ and 
\begin{equation}
\F(t) = F^{(n)} (\F(\tau_1), \dots, \F(\tau_n)),
\end{equation}
if $t = \Bplus(\tau_1,\dots,\tau_n)$.

B-series in any other pre-Lie algebra $(A, \tr)$ can be defined in the same way: by applying the unique pre-Lie algebra morphism $F:\PL \rightarrow A$ to the series (\ref{eq:plB}).
\begin{remark}\label{el.diff-grafting}
Since $\F: \PL \rightarrow \mathcal{X}(\RR^n)$ is a pre-Lie morphism, the trees associated to the derivatives of $y'(t) = F(y(t))$ can be generated by iterated grafting onto the one-vertex tree:
\[\text{The $n$ graftings}\quad  \ab \car (\ab \car (\ab \car \dots (\ab \car \ab) \dots )) \quad \text{corresponds to} \quad \frac{d^ny}{dt^n}.\]
This way of looking at elementary differentials will reappear in a different setting in Section \ref{sect:gengeomint}. 
\end{remark}

\begin{remark}\label{remark:magnus} The formulation of differential equations in terms of pre-Lie algebras has seen some use in numerical analysis. In \cite{ebrahimi-fard2009ama} Ebrahimi-Fard and Manchon rephrased differential equations of the type $X'(t) = A(t)X(t)$, where $X,A$ are linear operators in a vector space, as combinatorial equations in pre-Lie algebras. In this context they obtained an analogue of the Magnus expansion \cite{Magnus1954ote}, a series expansion of the solution to the equation  in the magma generated by monomials of pre-Lie elements. In this setting it becomes apparent that one can use the pre-Lie relation to cancel out some of the terms in the expansion, leading to a hitherto unknown reduction of the number of terms in the Magnus expansion.
\end{remark}

\section{Geometric numerical integration on manifolds}\label{sect:gengeomint}

Our objects of study are now dynamical systems evolving on \emph{manifolds}: 
\begin{equation}\label{manifolddiffeq}
y'=F(y), \hspace{0.6cm} y_0 \in M, \quad F \in \XM,
\end{equation}
where $M$ is a smooth manifold and $\XM$ denotes the vector fields on $M$. As in the previous chapter, the aim is to find good numerical approximations to the flow $\exp(tF) := \Psi_{t,F}$ of (\ref{manifolddiffeq}). The study of such systems comprises several different approaches: One simple way to attack the problem is to embed the manifold in $\mathbb{R}^N$, for some $N$, and use methods developed for $\mathbb{R}^N$ to solve the equation. But then the numerical flow of the method may drift off the manifold, and this can in some cases cause problems \cite{engo1998mas, iserles1997nmo, calvo1996rkm, iserles1996qna}. 

A more satisfying and often better way is to use methods that are intrinsic to the manifold, and not rely on any embedding, which is the approach taken by \emph{Lie group integrators} \cite{iserles2000lgm}. Consider for instance a system evolving on the manifold $S^3$. By embedding $S^3$ in $\mathbb{R}^4$ one can use numerical methods that approximate the flow of the system using the basic motions of translations in $\mathbb{R}^4$. Another approach is to use \emph{rotations} to move around $S^3$: $y_{n+1} = Q_n y_n$ where $Q_n$ are orthogonal matrices, i.e. to use the action of the Lie group $SO(4)$ on $S^3$. This illustrates the intrinsic approach, where we are guaranteed not to drift off $S^3$. Methods developed for manifolds include the Crouch--Grossman and RKMK-methods (and variants thereof) \cite{munthe-kaas1998rkm, munthe-kaas1999hor, crouch1993nio, owren1999rkm, engo2000otc}.

In this chapter we will study a generalization of B-series called \emph{Lie--Butcher series}. In analogy to the previous chapter we will look at the composition and substitution of Lie--Butcher series.

\subsection{Setting the stage: homogeneous manifolds and differential equations}\label{homogeneous}
The flows we would like to approximate evolve on smooth manifolds, and
so the tools of differential geometry play an important role. We will
not review the general theory of smooth manifolds here, but assume a
basic knowledge of differential geometry; for excellent introductions
see e.g. \cite{abraham1988mta, spivak2005aci, sharpe1997dg}. For a viewpoint oriented toward geometric numerical integration, see
\cite{iserles2000lgm}. More precisely, we will be working with smooth
manifolds equipped with transitive actions by Lie groups, so called
\emph{homogenous manifolds}, where the Lie group provides a way to
move around on the manifold.\footnote{Note that other manifolds with \emph{local} actions could also be considered, but to avoid unnecessary complications we only consider homogeneous manifolds.} Because the action is not in general free, the differential equation expressed on the Lie group is not in general unique.

\begin{definition}
An \emph{action of a Lie group} $G$ on a smooth manifold $M$ is a group homomorphism $\lambda: G \rightarrow \Diff(M)$, $g \mapsto \lambda_g$,  where $\Diff(M)$ is the group of diffeomorphisms on $M$. We will mostly write such an action as a map $\Lambda: G \times M \rightarrow M$.
\end{definition}
\noindent For convenience of notation we write $g$ for the diffeomorphism $\lambda_g$, and also $g \cdot m$ for $\lambda_g(m)$. The \emph{orbit} through a point $p \in M$ is the set $G\cdot p = \lambda_G(p)$. The action is called \emph{transitive} if the manifold $M$ is a single $G$-orbit. That is, if for all $p, q \in M$ there is a $g \in G$ so that $p = g \cdot q$. A manifold equipped with a transitive action by a Lie group $G$ is called a \emph{homogeneous manifold}. A consequence of this is that $M$ is diffeomorphic to the right cosets $G/G_x$ of $G$, where $G_x$ is the closed Lie subgroup of isotropies, $G_x = \{g \in G \,\,|\,\, gx = x\}$ (the point stabilizer): the smooth manifold structure of $G/G_x$ comes from the quotient map, and the diffeomorphism $F: G/G_x \rightarrow M$ is given by $F(gG_x) = g \cdot x$. The group $G_x$ is called {\it the} subgroup of isotropies because if $x'$ is another point in $G$, then $G_x$ and $G_{x'}$ are conjugate, and therefore isomorphic.

Some interesting examples of homogeneous manifolds are the spheres $S^n = SO(n+1)/SO(n)$. Other important examples come from Lie groups $G$ themselves. The action $\Lambda: G \times G \rightarrow G$ is then the Lie group multiplication. A (somewhat degenerate) example is the homogeneous manifold $(\RR^n, (\RR^n, +))$. Here the action of $\RR^n$ on itself is given by translations. The theory developed for homogeneous manifolds in this chapter will reduce to the theory developed in the previous chapter when applied to this particular case.

Actions by Lie groups on manifolds can be associated to actions by Lie algebras. Let $\Lambda: G\times M \rightarrow M$ be an action of $G$ on $M$. The associated Lie algebra action $\lambda_*: \g  \rightarrow \XM$ of $\g$ on $M$ is the Lie algebra anti-homomorphism defined by:
\begin{equation}\label{eq:infinitesimal}
\lambda_*(v)(p) = \left.\frac{d}{dt}\right|_{t=0} \Lambda(\exp(tV),p).
\end{equation}
This satisfies $\lambda_*([U,V]) = -[\lambda_*(U),\lambda_*(V)]$, where the bracket on the right is the Jacobi bracket of vector fields.
We sometimes write $v \cdot y$ for the element $\lambda_*(v)(y) \in T_yM$. The \emph{Lie--Palais theorem} \cite{palais1957agf} ensures us that as long as the Lie group $G$ is simply connected, then every action by $\g$ comes from an action by $G$.\footnote{If the Lie group is not simply connected, then we can only lift the $\g$-action to the universal covering group of $G$.} If $F \in \XM$ is a vector field, then an element $v$ so that $\lambda_*(v) = F$ is called an \emph{infinitesimal generator} for $F$.

\begin{remark}
In some cases it makes sense to use other maps $\phi: \g \rightarrow G$ (satisfying $\phi(0)=e$ and $\phi'(0)=V$) besides the exponential map to construct maps $\g \rightarrow \XM$ as in Equation (\ref{eq:infinitesimal}). An overview of various maps of this kind, and their usefulness, can be found in \cite{engo2000otc}. 
\end{remark}

\paragraph{Differential equations in homogeneous manifolds.} 
Consider the differential equation on a homogeneous manifold $(M,G,\lambda)$:
\begin{equation}\label{ivp2}
y' = F(y), \qquad y_0 \in M, \quad F: M \rightarrow TM.
\end{equation}
The solution is the flow $\Psi_{t, F} = \exp(tF)$ of the vector field $F$. 
The vector field can be written in terms of its infinitesimal generator as $F = \lambda_*(v): M \rightarrow TM$ for an element $v \in \g$, and the transitivity of the action also allows us to construct a map $f:M \rightarrow \g$ so that 
\begin{equation}
F(y) = \lambda_*(f(y))(y) = f(y)\cdot y
\end{equation}
Note that as long as the action is not free, this $f$ is not unique: if $f: M \rightarrow \g$ is such a map, then $f + i: M \rightarrow \g$, where $i(p)$ is in the isotropy subalgebra $\g_p$ of $\g$, is another map of the same type. This choice of isotropy class can be helpful when constructing numerical integrators \cite{lewis2002gia}. 

The differential equation (\ref{ivp2}) can be written as:
\begin{equation}\label{eq:diffeqLB}
y' = f(y) \cdot y, \quad \text{where} \quad f:M \rightarrow \g,
\end{equation}
and this is the type of differential equation we will consider in this chapter. Note that in the classical case of $(\RR^n, (\RR^n, +))$, the exponential is $\exp(v) = v$, and Equation \ref{eq:diffeqLB} reduces to the ordinary differential equation (\ref{ivp2}). We also note that the class contains the equations formulated in terms of \emph{frames}.

\begin{remark}[Frames and differential equations] In the literature for numerical integration of differential equations on manifolds the equations are often simplified by using a frame on the manifold \cite{owren1999rkm, owren2006ocf, celledoni2003cfl}.  A \emph{frame} is a set of vector fields $\{E_i\}$ that at each point on the manifold spans the tangent space at that point, so that any vector field $F$ can be written as $F=\sum_i f_i E_i$. The flow equation (\ref{ivp2}) for $F$ can then be written as 
\begin{equation}\label{eq:frameDE}
y' = \sum_i f_i(y) E_i(y), \quad \text{where} \quad f_i: M \rightarrow \RR \quad \text{are smooth}.
\end{equation}
If we write $\g \subset \XM$ for the Lie subalgebra generated by the frame vector fields $\{E_i\}$, and let $\lambda_*: \g \rightarrow \Diff(M)$ be as in (\ref{eq:infinitesimal}), we see that Equation (\ref{eq:frameDE}) is a special case of Equation (\ref{eq:diffeqLB}), with $f:M \rightarrow \g$ defined by $f(y) = \sum_i f_i(y) E_i$.
\end{remark}

\begin{remark}
In \cite{engo2000otc}, K. Eng\o{} discusses the general operation of `moving' differential equations between manifolds using equivariance of actions and relatedness of vector fields. In particular, every differential equation of the form (\ref{eq:diffeqLB}) is shown to be equivalent to a differential equation on $\g$. The following diagram from \cite{engo2000otc} summarizes this:
\begin{diagram}[labelstyle=\scriptstyle]
T\g &\rTo^{T(\exp)}& TG &\rTo^{T(\lambda_{\cdot}(p))}& TM \\
\uTo && \uTo && \uTo_{\lambda_*(v)(p)} \\
\g &\rTo_{\exp}& G &\rTo_{\lambda_{\cdot}(p)}& M
\end{diagram}
In other words, the differential equation on a homogeneous manifold $(M,G)$ is moved to the Lie group $G$ (the middle vertical arrow) and then to the Lie algebra $\g$ (the first vertical arrow). As before, the exponential map $\exp:\g \rightarrow G$ can in many cases be replaced by other maps. The construction of the vertical arrows can be found in \cite{engo2000otc}. This is the result exploited in the so-called RKMK methods \cite{munthe-kaas1995lbt, munthe-kaas1998rkm, munthe-kaas1999hor}.
\end{remark}

\subsection{Lie group integrators in applications} 
\begin{definition}
Given a smooth manifold $M$ and a Lie group $G$ with Lie algebra $\mathfrak{g}$ acting on $M$, consider a differential equation for $y(t)\in M$ written in terms of the infinitesimal action as
\begin{equation}\label{eq:lgeqn}
\dot{y}(t) = f(t,y)\cdot y, \quad y(0) = y_0,
\end{equation}
for a given function $f:\mathbb{R}\times M\rightarrow \mathfrak{g}$. Note that we now make the non-autonomousness explicit. A \emph{Lie group integrator} is a numerical time-stepping procedure for (\ref{eq:lgeqn}) based on intrinsic Lie group operations, such as exponentials, commutators and the group action on $M$. 
\end{definition}

Applications of LGI generally involve the following steps:
\begin{enumerate}
\item Choose a Lie group and Lie group action that can be computed fast and which captures some essential features of the problem to be solved. This is similar to the task of finding a preconditioner in iterative solution of linear algebraic equations.
\item Identify the Lie algebra, commutator and exponential map of the Lie group action. 
\item Write the differential equation in terms of the infinitesimal Lie algebra action, as in Equation~(\ref{eq:lgeqn}).
\item Choose a Lie group integrator, plug in all the building blocks, and solve the problem.
\end{enumerate}

\paragraph{Examples of Lie group integrators.} We list some important Lie group integrators:\\

\noindent{\bf Lie Euler: } \[\quad y_{n+1} = \exp(h f(t_n,y_n))\cdot y_n\]

\noindent{\bf Lie midpoint:} 
\begin{align*}K &= hf\left(t_n+ h/2,\exp\left(K/2\right)\cdot y_n\right)\\ y_{n+1} & =  \exp(K)\cdot y_n\end{align*}

\noindent{\bf Lie RK4:}
There are several similar ways of turning the classical RK4 method
into a fourth order Lie group integrator \cite{munthe-kaas1998rkm,munthe-kaas1999cia}. The following version requires only two commutators:
\begin{align*}
K_1 & =  hf(t_n,y_n)\\
K_2 & =  hf(t_n/2,\exp(K_1/2)\cdot y_n)\\
K_3 & =  hf\left(t_n+h/2, \exp(K_2/2-[K_1,K_2]/8)\cdot y_n\right)\\
K_4 & =  hf(t_n+h/2,\exp(K_3)\cdot y_n)\\
y_{n+1} & =  \exp\left(K_1/6+K_2/3+K_3/3+K_4/6 -[K_1,K_2]/3-[K_1,K_4]/12\right)\cdot y_n
\end{align*}

\noindent{\bf RKMK methods:} There is a general procedure to turn
any classical Runge--Kutta method into a Lie group integrator of the
same order \cite{munthe-kaas1999hor,iserles2000lgm}. 

Let $\{a_{i,j},b_j,c_i\}_{i,j=1}^s$ be the coefficients of a Runge--Kutta method of order $p$. The following method is a Lie
group method of order $p$.

\begin{proglist}
  for $i = 1,s$\+\\
     $U_i = \sum_{j=1}^s a_{i,j}{{{K}_j}}$\\
     $K_i = d\exp^{-1}_{U_i}\left(h f(\exp(U_i)\dpr y_n)\right)$
     %$\tilde{K}_i = d\exp^{-1}_{U_i}(K_i)$
\-\\
  end\\
  $y_{n+1} = {{\exp\left(\sum_{j=1}^s b_{j}{{{K}_j}}\right)\dpr y_n}}$ ,
\end{proglist}
where \[d\exp^{-1}_{U}(K)=\sum_{n=0}^\infty\frac{B_n}{n!} \ad^n(U)(K) =  { K - \frac12[U,K] + \frac{1}{12}[U,[U,K]] + \cdots}\]
is the inverse of the right trivialized tangent map of the exponential, see~\cite{iserles2000lgm}.
\ \\

\noindent{\bf Crouch--Grossman and commutator free methods:}  Commutators pose a  problem in the application of Lie group integrators to stiff problems, since the commutator often increases the stiffness of the equations dramatically. Crouch--Grossman~\cite{crouch1993nio,owren1999rkm}, and more generally commutator-free methods~\cite{celledoni2003cfl}, avoid them by doing basic time-stepping using a composition of exponentials. An example of such a method is CF4~\cite{celledoni2003cfl}:
\begin{align*}
K_1 & =  hf(t_n,y_n)\\
K_2 & =  hf(t_n/2,\exp(K_1/2)\cdot y_n)\\
K_3 & =  hf\left(t_n+h/2, \exp(K_2/2)\cdot y_n\right)\\
K_4 & =  hf(t_n+h/2,\exp(K_1/2)\cdot\exp(K_3-K_1/2)\cdot y_n)\\
y_{n+1} & =  \exp\left(K_1/4+K_2/6+K_3/6 -K_4/12\right)\cdot\\ & \hspace{0.5cm} \exp\left(K_2/6+K_3/6+K_4/4 -K_1/12\right)\cdot y_n
\end{align*}

\noindent{\bf Magnus methods:} In the case where $f(t,y)=f(t)$ is a
function of time alone, then~(\ref{eq:lgeqn}) is called an equation of
\emph{Lie type}. Specialized numerical methods have been developed for
such problems~\cite{iserles1999sld, iserles1999imm, blanes2009tme}. Explicit Magnus methods can achieve order 2p using only p function evaluations, and they are also easily designed to be time symmetric.
%\end{example}

\paragraph{Examples of group actions.} Some group actions of interest
when applying Lie group integrators:\\

\noindent{\bf Rotational problems:} Consider a differential equation $\dot{y}(t) = v(t)\times y(t)$, where $y,v\in \mathbb{R}^3$ and $||y(0)||=1$. Since $||y(t)||=1$ for all $t$, we can take $M$ to be the surface of the unit sphere. Let $G=SO(3)$ be the special orthogonal group, consisting of all orthogonal matrices with determinant 1. Let $\gamma(t)\in G$ be a curve such that $\gamma(0)=e$. By differentiating $\gamma(t)^T\gamma(t)=e$, we find that $\dot{\gamma}(0)^T+\dot{\gamma}(0)=0$, thus $\mathfrak{g}=\mathfrak{so}(3)$, the set of all skew-symmetric $3\times 3$ matrices. The infinitesimal Lie algebra action is left multiplication with a skew matrix, the commutator is the matrix commutator, and the exponential map is the matrix exponential. Written in terms of the infinitesimal Lie algebra action, the differential equation becomes $\dot{y}(t) = \widehat{v(y)} y$, and we may apply any Lie group integrator. Note that for low dimensional rotational problems, all basic operations can be computed fast using Rodrigues type formulas~\cite{iserles2000lgm}.\\

\noindent{\bf Isospectral action:} Isospectral differential equations are matrix valued equations where the eigenvalues are first integrals (invariants of motion). Consider $M= \mathbb{R}^{n\times n}$ and the action of $G=SO(n)$ on $M$ by similarity transforms, i.e.\ for $a\in G$ and $y\in M$ we define $a\cdot y = aya^T$. By differentiation of the action we find the infinitesimal action for $V\in\mathfrak{g}=\mathfrak{so}(n)$ as $V\cdot y = Vy - yV$, thus for this action~(\ref{eq:lgeqn}) becomes \[\dot{y}(t) = f(t,y)\cdot y = f(t,y)y-yf(t,y),\] where $f\colon \mathbb{R}\times M\rightarrow \mathfrak{g}$. See~\cite{calvo1997nso,iserles2000lgm} for more details.\\

\noindent{\bf Affine action:} Let $G=Gl(n)\rtimes \mathbb{R}^n$ be the \emph{affine linear group}, consisting of all pairs $a,b$ where $a\in \mathbb{R}^{n\times n}$ is an invertible matrix and $b\in \mathbb{R}^n$ is a vector. The \emph{affine action} of $G$ on $M=\mathbb{R}^n$ is $(a,b)\cdot y = ay+b$.
The Lie algebra of $G$ is $\mathfrak{g} = \mathfrak{gl}(n)\rtimes \mathbb{R}^n$, i.e.\ $\mathfrak{g}$ consists of all pairs $(V,b)$ where $V\in \mathbb{R}^{n\times n}$ and $b\in \mathbb{R}^n$. The infinitesimal action is given as $(V,b)\cdot y = Vy + b$. This action is useful for differential equations of the form $\dot{y}(t) = L(t)y + N(y)$, where $L(t)$ is a stiff linear part and $N$ is a non-stiff non-linear part. Such equations are cast in the form~(\ref{eq:lgeqn}) by choosing $f(t,y) = (L(t),N(y))$. Applications of Lie group integrators to such problems is closely related to exponential integrators. In this case it is important to use a commutator-free Lie group method.
\\

\noindent{\bf Coadjoint action:} 
Many problems of computational mechanics are naturally formulated as
\emph{Lie--Poisson systems}, evolving on coadjoint orbits of the dual of a Lie
algebra~\cite{marsden1999itm}. Lie group integrators based on the
coadjoint action of a Lie group on the dual of its Lie algebra are
discussed in \cite{engo2002nio}.
\\

\noindent{\bf Classical integrators as Lie group integrators:} The simplest of all group actions in our setting arises when $G=M=\mathbb{R}^n$. We can then use vector addition as group operation and group action. From the definitions we find that in this case $\mathfrak{g} = \mathbb{R}^n$, the commutator is 0, and the exponential map is the identity map from $\mathbb{R}^n$ to itself. The infinitesimal Lie algebra action becomes $V\cdot y = V$, thus~(\ref{eq:lgeqn}) reduces to $\dot{y}(t) = f(t,y)$, where $f(t,y)\in \mathbb{R}^n$. We see that classical integration methods are special cases of Lie group integrators, and all the examples of methods above reduce to well-known Runge--Kutta methods.
%\end{example}

\subsection{Geometry meets algebra}\label{sect:geometry}
We will discuss how the B-series of Section~\ref{sect:classgeomint} can be generalized to \emph{Lie--Butcher series} for analyzing exact and approximate flows on manifolds.
The basic building blocks of the numerical methods discussed above are commutators in the Lie algebra of frame vector fields $\g$, flows of frozen vector fields, evaluation of frozen vector fields, and parallel transport of tangent vectors. The Lie algebra $\g$ defines an \emph{absolute parallelism} on the Lie group, which yields a parallel transport of  tangent vectors \cite{abraham1988mta}. For a vector field $F: \M\rightarrow T\M$ represented by a function $f: \M\rightarrow \g$, parallel transport between $T_{gy}\M$ and $T_y\M$ is defined as $\tau_g f(y)= f(gy)$. This transport is independent of the path between $gy$ and $y$, and hence it is a parallel transport induced by a \emph{flat connection}. Furthermore, we will see that this connection has \emph{constant torsion}. This is the geometric setting of Lie--Butcher series: a manifold with a connection which is flat with constant torsion, giving rise to a \emph{post--Lie algebra}~\cite{vallette2007hog,munthe-kaas2012opl} of vector fields on $\M$. The enveloping algebra of a post--Lie algebra is called a \emph{D-algebra} \cite{munthe-kaas2008oth,lundervold2009hao,lundervoldthesis}.

Let $(G,\M)$ be a homogeneous manifold and let $\g$ denote the Lie algebra of $G$, which can be identified with right invariant vector fields (i.e.\ invariant derivations on $C^\infty(\M)$) via~(\ref{eq:infinitesimal}). We let $U(\g)$ denote the universal enveloping algebra of $\g$, which similarly can be identified with higher order right invariant derivations on $C^\infty(\M)$. Let $\one$ denote the identity operator on $C^\infty(\M)$.
We define
\begin{align}
\begin{split}
\UgM & :=  C^\infty(M)\tpr_\RR U(\g)\\
\gM & :=  C^\infty(M)\tpr_\RR \g \subset \UgM.
\end{split}
\end{align}
Let $\left\{\partial_i\right\}_i$ be a basis for $\g$. Then $f\in \gM$ can be written as $f = \sum_i f^i\tpr \partial_i$. This represents a function $f\colon \M\rightarrow\g$ as $f(x) = \sum_i f^i(x)\partial_i\in \g$, and $f$ also acts as a derivation
$f[g]:=\sum_i f^i\partial_i g$ for $g\in C^\infty(\M)$. Similarly,
higher order derivations on $C^\infty(\M)$ can be represented by
\begin{equation}
h =\sum_I h^I \tpr\partial_I =
\sum_{i,j,k\ldots}h^{ijk\ldots}\tpr \partial_i\partial_j\partial_k\ldots\in
\UgM,
\end{equation}
where $I=(i,j,k,\ldots)$ is a multi-index.

The space $\UgM$ is equipped with two operations: \emph{frozen composition} $g,h\mapsto gh$ and \emph{covariant derivation} $g,h\mapsto g[h]$ defined 
for $g,h\in \UgM$,  $g=\sum_I g^I\tpr\partial_I$, $h=\sum_J h^J\tpr\partial_J$ as
\begin{align}
gh & =  \sum_{I,J}g^Ih^J\tpr \partial_I\partial_J\\
g[h] & =  \sum_{I,J} g^I\partial_Ih^J\tpr \partial_J.
\end{align}
Note that these operations are independent of the choice of basis for $\g$.
The covariant derivation $g[h]$ can be understood as a (higher order) derivative of $h$ as it moves under parallel transport defined by the absolute parallelism.
We will frequently apply the alternative notation 
\begin{equation}
g\tr h := g[h]
\end{equation}
to emphasize the similarity between this operation and the product $\tr$ in a Pre-Lie algebra.

Let $\D(\UgM)\subset\UgM$ denote the (first order) \emph{derivations}, defined as
 \[\D(\UgM) = \stset{f\in\UgM}{f[gh] = (f[g])h + g(f[h])\quad\mbox{for all $g,h \in \UgM$}}.\]
It can be shown that
\begin{equation}\D(\UgM)=\gM .
\end{equation}
$\UgM$ with the operations of frozen composition and covariant derivation satisfies the following fundamental relationships:
For any derivation $f \in \D(\UgM)$ and any $g,h\in \UgM$ we have
\begin{align}
\begin{split}
 g[f] &\in \mathcal{D}(\UgM)\\
 f[g[h]] &= (fg)[h] + (f[g])[h].
\end{split}
\end{align}
Such an algebraic structure is called a \emph{D-algebra} in
\cite{munthe-kaas2008oth}. We shall see that the \emph{free D-algebra} is the algebra of forests of ordered trees, where $f,g\mapsto fg$ is the concatenation of forests and $f,g\mapsto f[g]$ is left grafting. 

\subsection{Ordered trees and D-algebras.}
The set
\[\OT = \{\ab,\aabb,\aababb,\aaabbb, \aabababb, \aabaabbb, \aaabbabb, \aaababbb,\ldots \}.\]
of ordered rooted trees consists of {all} \emph{planar} rooted trees. In other words, an ordered rooted tree is a tree $\tau$ together with a chosen order of the branches connected to each vertex of $\tau$. Unlike the set $\T \subset \OT$ of rooted trees, we do not identify trees who differ in the order of their branches. 
%We%write $\ROT$ for the free $\RR$-vector space over the set $\OT$,
%i.e. all $\RR$-linear combinations of ordered trees. 

The set of ordered words of elements from $\OT$ is denoted by $\OF$, and also includes the empty word. $\OF$ is called the set of \emph{ordered forests}, and we write $\OF_{\C}$ for the set of forests colored by $\C$. Let $\N = \RR\langle \OT\rangle$ be the non-commutative polynomials over
$\OT$. The linear dual $\N^* := \Hom(\N, \RR)$ is identified with the
infinite combinations of forests, and we write $\langle \cdot, \cdot
\rangle$ for the pairing making forests orthogonal. That is,
$\langle \omega_1, \omega_2\rangle = \delta_{\omega_1, \omega_2}$, for
all $\omega_1, \omega_2 \in \OF$.

It is sometimes convenient to allow the trees to be \emph{decorated} by a set $\C$, often called the set of colors. This is done via a map from the vertices of the tree to the set $\C$. We write $\OT_{\C}$ and $\OF_{\C}$ for the set of trees and forests colored by $\C$.

A basic operation on $\N$ is the \emph{left grafting product} $\car: \N \otimes \N \rightarrow \N$ of \cite{munthe-kaas2008oth}. It is defined recursively by
\begin{align}
\begin{split}
&\one\car\omega = \omega \\
&\omega\car\one = 0 \\
&\omega\car\ab  =  \Bplus(\omega),\\
&\tau\car\omega_1 \omega_2 = (\tau\car\omega_1)\omega_2 + \omega_1(\tau\car\omega_2) \\
&(\tau\omega)\car\omega_1 =
\tau\car(\omega\car\omega_1) - (\tau\car\omega)\car\omega_1,
\end{split}
\end{align}
where $\tau$ is a tree and $\omega_1$, $\omega_2$ are
forests. 

The grafting product can be used to define the \emph{Grossman-Larson
  product} (GL product) $\bpr: \N \otimes \N \rightarrow \N$: 
\begin{align}
\Bplus(\omega_1 \bpr \omega_2) = \omega_1\car \Bplus(\omega_2),
\end{align}
extended by linearity.

Concatenation and left grafting gives $\N$ the structure of a
\emph{D-algebra}, as defined in \cite{munthe-kaas2008oth} (see also
\cite{lundervold2009hao, lundervold2010bea, munthe-kaas2012opl}), where
the composition $\tr$ is left grafting.

\begin{definition}\label{Dalg}
Let $A$ be a unital associative algebra with product $f,g \mapsto fg$
and unit $\one$, equipped with a non-associative composition $\tr: A
\otimes A \rightarrow A$ such that $\one\tr g = g$ for all $g\in
A$. Write $\mathcal{D}(A)$ for the set of derivations:
$$\mathcal{D}(A) = \{f\in A \,\,|\,\, f\tr (gh) = (f\tr g)h + g(f\tr h) \hspace{0.2cm} \text{for all } g,h \in A\}.$$ Then $A$ is called a \emph{D-algebra} if for any derivation $f \in \mathcal{D}(A)$ and any $g,h\in A$ we have
\begin{align*}
&\text{(i)} g\tr f \in \mathcal{D}(A)\\
&\text{(ii)} f\tr(g\tr h) = (fg)\tr h + (f\tr g)\tr h.
\end{align*}
\end{definition}
In \cite{munthe-kaas2008oth} it was shown that the D-algebra $\N$ is the \emph{free} D-algebra:
\begin{theorem}[\cite{munthe-kaas2008oth}]\label{universalD-alg} The vector space $\N = \fieldk\langle\OT_{\C}\rangle$ is the free D-algebra over the set $\C$. That is, for any $D$-algebra $\A$ and any map $\nu\colon\C\rightarrow D(\A)$ there exists a unique D-algebra homomorphism $\F_\nu\colon N\rightarrow \A$ such that $\F_\nu(c)  = \nu(c)$ for all $c\in \C$.
\begin{diagram}[labelstyle=\scriptstyle]
\C &\rInto& \N \\
\dTo^{\nu} && \dTo_{\exists\,! \,\, \F_{\nu}} \\
D(\A) &\rInto& \A
\end{diagram}
\end{theorem}
A D-algebra homomorphism between two D-algebras $A$ and $B$ is an algebra morphism $F:A \rightarrow B$ such that $F(\mathcal{D}(A)) \subset \mathcal{D}(B)$, and $F(a[b]) = F(a)[F(b)]$. 

By applying this theorem to the D-algebra $\UgM$ of differential operators on
a homogeneous manifold (defined in Section \ref{sect:geometry}), we will construct elementary differentials and the Lie--Butcher series (Definition \ref{LBelementdiff} and Definition \ref{LBseries}). 

%This theorem enables us to define elementary differentials and Lie--Butcher series by applying it to the case where $\A$ is the D-algebra $U(\g)$ of differential operators. Recall that a vector field (or, in other words, a first-order differential operator) $F$ on a homogeneous manifold $(M,G)$ can be represented as a function $f:M \rightarrow \g$. Similarly, all higher order differential operators on $M$ can be represented as functions from $M$ to the universal enveloping algebra $U(\g)$ of $\g$.

%\begin{theorem}[\cite{munthe-kaas2008oth}]Let $(M,G)$ be a homogeneous manifold and let $\g$ denote the Lie algebra of $G$. Let $U(\g)$ denote the universal enveloping algebra of $\g$, consisting of all higher order differential operators on $M$, and extend its structure to $C^{\infty}(M,U(\g)) =: U(\g)^M$ via 
%\begin{equation}
%F[G](p) := (F(p)[G])(p), \quad FG(p) := F(p)G(p).
%\end{equation}
%These two operations give $U(\g)^M$ the structure of a D-algebra.
%\end{theorem}

\paragraph{Post-Lie algebras.} In Section \ref{sect:geometry} we noted
that the derivations in the D-algebra $\UgM$ of differential operators
could be identified with $\gM$. In general, the derivations in a
D-algebra form what is called a \emph{post-Lie algebra}, and the
D-algebra can be identified with the universal enveloping algebra of
its post-Lie algebra of derivations. This point of view is developed in
\cite{munthe-kaas2012opl}, and is also being studied further in an
ongoing project \cite{ebrahimi-fard2011otp} where the \emph{operad}
behind post-Lie and D-algebras (also called \emph{post-associative
  algebras}) is explored. Post-Lie algebras were introduced
independently by Vallette in \cite{vallette2007hog}, in a different context.

%In \cite{munthe-kaas2012opl} the authors developed a more refined view of D-algebras, where the D-algebras are enveloping algebras of \emph{post-Lie algebras} (post-Lie algebras were also introduced independently by Vallette in \cite{vallette2007hog}). This point of view is currently being studied further in an ongoing project \cite{ebrahimi-fard2011otp}, where the \emph{operad} behind post-Lie and D-algebras (also called \emph{post-associative algebras}) is explored.

\begin{definition} A \emph{post-Lie algebra} is a Lie algebra $(A, [\cdot,\cdot])$ equipped with a non-commutative, non-associative product $\triangleright: A \otimes A \rightarrow A$ satisfying:
\begin{align}
x \tr [y,z] &= [x\tr y, z] + [y, x \tr z] \hspace{1cm} \text{(derivation property)}\\
[x,y] \tr z &= a_{\tr}(x,y,z) - a_{\tr}(y,x,z),\label{pL:assoc}
\end{align}
where $a_{\tr}(x,y,z)$ is the associator $a_{\tr}(x,y,z) = x \tr (y \tr z) - (x\tr y) \tr z$.
\end{definition}
\noindent Notice that relation (\ref{pL:assoc}) implies that a pre-Lie algebra (Section \ref{sect:prelie}) is a post-Lie algebra with vanishing bracket. In \cite{munthe-kaas2012opl} it is shown that the free Lie algebra over rooted trees colored by a set $\C$, equipped with a post-Lie operation derived from grafting of trees, is the \emph{free post-Lie algebra}, and that its universal enveloping algebra is the free D-algebra defined above. We will not pursue this point of view in our present study of Lie--Butcher series.

\subsection{Lie--Butcher series}

Analogous to the B-series of Section \ref{sect:classgeomint}, the Lie--Butcher series can be used to represent flows -- numerical or exact -- on homogeneous manifolds. To achieve this one combines the concept of \emph{Lie series} in free Lie algebras with ideas from the theory of B-series. An exposition of free Lie algebras and Lie series can be found in the book \cite{reutenauer93fla} by Reutenauer. 

The \emph{free Lie algebra} $\FLA(A)$ over a set $A$ of generators is the closure of the generators under commutation and linear combination. In particular, we have the free Lie algebra $\FLA(\OT)$ over the set of ordered rooted trees. A \emph{Lie series} is a series expansion:
\begin{equation}
S = \sum_{n \geq 0} S_n,
\end{equation}
where each homogeneous component is an element of $\FLA(\OT)$, i.e. the $S_n$'s are \emph{Lie polynomials}.

A Lie series of particular interest to us appears when computing the pullback of functions along flows of vector fields on homogeneous manifolds. Let $F \in \XM$ be a vector field with flow $\Phi_{t,F}$, and $\psi:M \rightarrow \g$ a function. Then
\begin{equation}
\left.\frac{d}{dt}\right|_{t=0}\Phi_{t,F}^*\psi = F[\psi].
\end{equation}
The Taylor expansion of $\Phi_{t,F}^*\psi$ around $0$ therefore takes the form of a Lie series
\begin{align}\label{LieSeries}
\begin{split}
\Phi_{t,F}^*\psi &= \sum_{n=0}^{\infty}
\frac{t^n}{n!}\left(\left.\frac{\partial^n}{\partial 
t^n}\right|_{t=0} \Phi_{t,F}^*\psi \right)\\
&= \psi + tF[\psi]  + \frac{t^2}{2!} F[F[\psi]] + \frac{t^3}{3!}F[F[F[\psi]]] + \cdots.
\end{split}
\end{align}

\paragraph{Bell polynomials.} The higher order derivatives of the
pullbacks can be written in terms of non-commutative analogs of the
classical Bell polynomials of \cite{bell1927pp}. These polynomials
have also appeared in \cite{munthe-kaas1995lbt, schimming1996nbp, munthe-kaas1998rkm, lundervold2009hao}.

\begin{definition}
Let $D = \RR\langle \I \rangle$ be the free associative algebra over an alphabet $\I = \{d_i \}$, and let $\partial: D \rightarrow D$ denote the derivation given by $\partial(d_i) = d_{i+1}$. The \emph{non-commutative Bell polynomials} $B_n = B_n(d_1,\dots, d_n) \in \RR\langle \I\rangle$ are defined by the recursion
\begin{align}
\begin{split}
B_0 &= \one \\
B_n &= (d_1 + \partial)B_{n-1}, \quad n > 0.
\end{split}
\end{align} 
\end{definition}
The first few are:
 \begin{align*}
B_0 & =  \one\\
B_1 & =  d_1\\
B_2 & =  d_1^2 + d_2 \\
B_3 & =  d_1^3 + 2d_1d_2 + d_2d_1+d_3\\
B_4 & =  d_1^4+ 3d_1^2 d_2 +2d_1d_2d_1 + d_2d_1^2 +3d_1d_3+  d_3d_1 + 3d_2^2 + d_4.
\end{align*}

\begin{theorem}[{\cite{munthe-kaas1995lbt, lundervold2009hao}}] The derivatives of the pullback of a function $\psi$ along the time-dependent flow $\Phi_{t,F}$ can be written as:
\begin{equation}\label{eq:Bnderiv1} \frac{d^n}{d t^n}\Phi_{t,F}^*\psi = B_n(F)[\psi] ,\end{equation}
where $B_n(F_t)$ is the image of the Bell polynomials $B_n$ under the homomorphism given by $d_i\mapsto F^{(i-1)}$ ($(i-1)$th derivative). In particular
\begin{equation}
\left.\label{eq:Bnderiv} \frac{d^n}{d t^n}\right|_{t=0}\Phi_{t,F_t}^*\psi = B_n(F_1,\ldots,F_n)[\psi] =: B_n(F_i)[\psi],
\end{equation}
where $F_{n+1} = d^n/dt^n|_{t=0} F$.
\end{theorem}
This result allows us to obtain a  Lie series corresponding to (\ref{LieSeries}) for the case when $F$ is non-autonomous \cite{munthe-kaas1995lbt}:
\begin{equation}
\Phi^*_{t,F} \psi % = \sum_{n=0}^{\infty} F^n[\psi] \frac{t^n}{n!} [GJELDER KUN OM F ER AUTONOM]
= \sum_{n=0}^{\infty} B_n(F_i)[\psi] \frac{t^n}{n!}.
\end{equation}
%where $F^n$ iterated application of $F$, as in Equation (\ref{LieSeries}). 

\begin{remark} It is well known that the classical Bell polynomials \cite{bell1927pp} can be defined in terms of determinants. As an interesting side note, the non-commutative Bell polynomials can be defined in the same way, only now in terms of a non-commutative analog of the determinant: the \emph{quasideterminants} of Gelfand and Retakh (\cite{gelfand1991dom}, see also \cite{gelfand2005q}). For example, we have  

\begin{align*}
\det\begin{bmatrix}
x_1 & -1 & 0  \\  \\
{3-1 \choose 1} x_2 & x_1 & -1  \\  \\
\text{\circled{${3-1 \choose 2}x_3$}} & {3-2 \choose 1} x_2  & x_1 
\end{bmatrix}
&=
\det\begin{bmatrix}
x_1 & -1 & 0  \\  \\
2 x_2 & x_1 & -1  \\  \\
\text{\circled{$x_3$}} &  x_2  & x_1 
\end{bmatrix}\\
&\\
&= x_1^3+2x_1 x_2 + x_2 x_1 + x_3\\
&= B_3,
\end{align*}
where $\det$ denotes the quasideterminant, computed at the circled element. See \cite{gelfand2005q} for details about the computation and properties of quasideterminants. \end{remark}

The non-commutative \emph{partial Bell polynomials} $B_{n,k}:= B_{n,k}(d_1,\ldots,d_{n-k+1})$ are defined as the part of $B_n$ consisting of words $\omega$ of length $k>0$, 
e.g. $B_{4,3} = 3d_1^2 d_2 +2d_1d_2d_1 + d_2d_1^2$. Thus
\begin{equation}
B_n = \sum_{k=1}^n  B_{n,k} .
\end{equation}

\paragraph{A Fa\`a di Bruno bialgebra.} The non-commutative \emph{Dynkin--Fa\`a di Bruno bialgebra} $\Hfdb$ is obtained by using the algebra structure of $\Hfdb$ and defining the coproduct $\cpfdb$ as
\begin{align}\label{eq:fdb2}
\begin{split}
\cpfdb(\one)  & =  \one\tpr \one\\
\cpfdb(d_n) & =  \sum_{k=1}^n B_{n,k}\tpr d_k .
\end{split}
\end{align}
This extends to all of $\Hfdb$ by the product rule $\cpfdb(d_i d_j) =
\cpfdb(d_i)\cpfdb(d_j)$. For example
\begin{align*}
\cpfdb(d_1) & =  d_1\tpr d_1\\ 
\cpfdb(d_2) & =  d_1^2\tpr d_2 + d_2\tpr d_1\\
\cpfdb(d_1d_2) & =  d_1^3\tpr d_1d_2 + d_1d_2\tpr d_1^2.
\end{align*}
Note that the coproduct is not graded by $|\cdot|$

\begin{lemma}[{\cite{lundervold2009hao}}]~\label{coprod of Bell} The coproduct of the partial Bell polynomials is:
\begin{equation}\label{eq:fdb3}
\cpfdb(B_{n,k})  =  \sum_{\ell=1}^n B_{n,\ell}\tpr B_{\ell,k} .
\end{equation}
\end{lemma}
\noindent Note that  $B_{n,1}=d_n$, so (\ref{eq:fdb2}) is a special case of (\ref{eq:fdb3}). Summing the partial $B_{n,k}$ over $k$, we find the coproduct of the full Bell polynomials:
\begin{equation}
\cpfdb(B_{n})  =  \sum_{k=1}^n B_{n,k}\tpr B_{k}.
\end{equation}
\noindent Using Lemma \ref{coprod of Bell} and the fact that $B_{n,k} = 0$ for $k > n$, one can  show that $\Hfdb$ is a bialgebra.

\begin{proposition}[{\cite{lundervold2009hao}}] $\Hfdb = \RR\langle \I \rangle$ with the non-commutative concatenation product and the coproduct $\cpfdb$ form a bialgebra $\Hfdb$, which is neither commutative nor cocommutative. 
\end{proposition}

\paragraph{Lie--Butcher series.}
The Lie--series (\ref{LieSeries}) can also be written as the \emph{Lie--Butcher series} for the exact flow. In general, the Lie--Butcher series $\Bs_f(\alpha)$ are constructed to represent flows given by $y_0 \mapsto y_t = \Psi_{t}(y_0)$: 
\begin{equation}
\Psi_{t}(y(t)) = \Bs_f(\alpha)[\Psi_t](y_0).
\end{equation}
Before giving the definition of Lie--Butcher series we define the elementary differentials of a vector field $F$:

\begin{definition}\label{LBelementdiff}
Let $\F_f:\N \rightarrow U(\g)^M$ be the unique D-algebra morphism given by Theorem \ref{universalD-alg} by associating \,$\ab$\, to a vector field $f:M \rightarrow \g$. This is called the \emph{elementary differentials} of the vector field $f$. 
\end{definition}
\noindent Note that $\F_f:\N \rightarrow U(\g)^M$ is given recursively by
\begin{itemize}
\item[(i)] $\F_{f}(\one) = \one$
\item[(ii)] $\F_{f}(\Bplus(\omega)) = \F_{f}(\omega)[f]$
\item[(iii)] $\F_{f}(\omega_1\omega_2) = \F_{f}(\omega_1)\F_{f}(\omega_2)$
\end{itemize}
\noindent The general Lie--Butcher series are expansions of elementary differentials indexed over ordered rooted forests.
\begin{definition}\label{LBseries}
A \emph{Lie--Butcher series} (LB-series) is a formal series expansion over $U(\g)^M$:
\begin{equation}
\Bs_f(\alpha) = \sum_{\omega \in \OF} h^{|\omega|}\alpha(\omega)\F_f(\omega),
\end{equation}
where $\alpha: \N \rightarrow \RR$.
\end{definition}
\noindent It turns out \cite{lundervold2009hao} that the Lie series (\ref{LieSeries}) can be written as 
\begin{equation}
\Phi^*_{t,f} \psi = \sum_{\omega\in\OT} \gamma(\omega) \F_{f}(\omega),
\end{equation}
where $\gamma$ are the coefficients appearing when iteratively (left) grafting $\ab$ onto $\ab$.
This is the Lie--Butcher series for the exact flow. 

To understand how Lie--Butcher series can be used to represent
numerical flows we conduct a closer study of the coefficients $\alpha:
\N \rightarrow \RR$, and understand them as characters in a certain
Hopf algebra. This Hopf algebra allows us to formulate the concept of
\emph{composition} of LB-series.

\subsection{Composition of Lie--Butcher series}\label{sect:compLB}
We would like to understand the result of composing LB-series methods in a similar way as we did for B-series methods in Section \ref{GI:composition}. The basic problem is to determine whether the method  $\Phi$ resulting from composing two methods $\Phi^2 \circ \Phi^1$--both given by LB-series--is another LB-series, and in that case, what its coefficients are. Just as there is a Hopf algebra governing composition of B-series (the BCK Hopf algebra discussed in Section \ref{GI:composition}), there is a Hopf algebra $\HMKW$ behind the composition of LB-series. This Hopf algebra was  studied in \cite{munthe-kaas2008oth}, where its properties and its relation to the BCK Hopf algebra was explored. An introduction can also be found in \cite{lundervold2009hao}. This 
Hopf algebra is the dual of a version of Grossman and Larsons Hopf algebras in the case of ordered trees.

\paragraph{The Hopf algebra of composition.} As a vector space $\HMKW$ is equal to $\N$: $\HMKW = \RR\langle \OT\rangle$. The product is given by \emph{shuffling}: 
\begin{align}
\begin{split}
&\one \shuffle \omega = \omega = \omega \shuffle \one\\
&(\tau_1\omega_1) \shuffle (\tau_2 \omega_2) = \tau_1(\omega_1 \shuffle \tau_2 \omega_2) + \tau_2(\tau_1\omega_1 \shuffle \omega_2)
\end{split}
\end{align}
where $\tau_1, \tau_2 \in \OT$ and $\omega_1, \omega_2 \in \OF$. The coproduct is given recursively by $\DeltaMKW(\one) = \one \otimes \one$ and 
\begin{equation}
\DeltaMKW(\omega \tau) = \omega \tau \otimes \one + \DeltaMKW(\omega) \shuffle \, \cdot \, (I \otimes B^+) \DeltaMKW(B^-(\tau)),
\end{equation}
 where $\tau \in \OT$, $\omega \in \OF$. Here $\shuffle \,\cdot: \N^{\otimes 4} \rightarrow \N \otimes \N$ denotes shuffle on the left and concatenation on the right: $(\omega_1 \otimes \omega_2) \shuffle \, \cdot\, (\omega_3 \otimes \omega_4) = (\omega_1 \shuffle \omega_3) \otimes (\omega_2 \omega_4).$

Note that the shuffle product also gives rise to the
\emph{shuffle Hopf algebra} $\Hsh$, whose coproduct is given by
deconcatenation \cite{reutenauer93fla}:
\begin{equation}
\DeltaC(\omega) = \one \otimes \omega + \omega \otimes \one + \sum_{i=1}^{n-1} \tau_1 \cdots \tau_i \otimes \tau_ {i+1} \cdots \tau_n,
\end{equation}
where $\omega = \tau_1 \cdots \tau_n$.

The set of ordered forests can be generated recursively from the empty forest $\one$ by a \emph{magmatic} operation $\times: \N \times \N \rightarrow \N$ on $\N$, given by $\mu_{\times}(\omega_1,\omega_2) = \omega_1 \Bplus(\omega_2)$. For a forest $\omega$, write $\omega_L$ and $\omega_R$ for the \emph{left-} and \emph{right part}: $\omega = \omega_L \times \omega_R$. The above operations can be written recursively in terms of this operation:
\begin{itemize}
\item[] Concatenation: $\omega\,\one = \one\,\omega = \omega$, and $(\omega_1\times \omega_2)\,\omega_3 = \omega_1 \times (\omega_2 \,\omega_3)$.
\item[]Shuffle: $\omega \shuffle \one = \one \shuffle \omega = \omega$, and $\omega_1 \shuffle \omega_2 = (\omega_1 \shuffle \omega_{2L}) \times \omega_{2R} + (\omega_{1L} \shuffle \omega_2) \times \omega_{1R}$
\item[]Coproduct $\DeltaMKW(\one) = \one \otimes \one$, and $\DeltaMKW(\omega) = \omega \otimes \one + \DeltaMKW(\omega_L) \shuffle \times \DeltaMKW(\omega_R)$
\end{itemize}

The coproduct can also be written in terms of \emph{left admissible cuts}, analogous to the coproduct in $\BCK$ (Theorem \ref{thm:BCKcoprod}):

\begin{theorem}[{\cite{munthe-kaas2008oth}}]
The coproduct in $\HMKW$ can be written as
\begin{equation}\label{MKWcoprod_cuts}
\DeltaMKW(\omega) = \sum_{c \in \FLAC(\omega)} P^c(\omega) \tpr R^c(\omega),
\end{equation}
where $\omega$ is a forest in $\OT$. Here $\FLAC(\omega)$ consists of all left admissible cuts of $\omega$, including the full cut.
\end{theorem}
A \emph{left admissible cut} differs from the admissible cuts defined in Section \ref{GI:composition}: an \emph{elementary cut} $c$ of a tree $\tau$ is a selection of edges to be removed from $\tau$, chosen in such a way that if an edge $e$ is removed, then all the branches on the same level and to the left of $e$ must also be removed. A cut results in a collection of trees concatenated together to form a forest $P^{c}_{el}(\tau)$ (the \emph{pruned part}), and a remaining tree $R^c_{el}(\tau)$, containing the root. A left admissible cut $c = \{c_1, \dots, c_n\}$ on $\tau$ is a collection of such elementary cuts, with the property that any path from the root to any vertex crosses at most one cut $c_i$. The pruned parts from each cut together form the pruned part $P^c(\tau)$ of the left admissible cut, where the parts coming from different cuts are shuffled together. We also include the \emph{full cut} and the \emph{empty cut}, which results in $P^c(\tau) = \tau$ and $P^c(\tau)=\one$, respectively. The cutting operation is extended to forests $\omega$ as follows: apply the $\Bplus$ operation to $\omega$ to get a tree, cut this according to the above rules, without using the cut removing all the edges coming out of the root, and, finally, remove the added root from $R^c(\omega)$.

See Table \ref{MKWcoprod} for some examples of the coproduct $\DeltaMKW$, and see \cite{munthe-kaas2008oth} or \cite{lundervold2009hao} for further examples and properties of $\HMKW$.
\begin{table}[h! ]
  \centering
  \begin{equation*}
    \begin{array}{c@{\,\,}|@{\quad}l}
      \hline \\[-2mm]
      \omega & \DeltaMKW(\omega) \\[1mm]
      \hline \\[-2mm]
      \one & \one\otimes \one \\[1mm]
      \ab & \ab\otimes \one + \one \otimes \ab \\[2mm]
      \ab\,\ab & \ab\,\ab\otimes \one + \ab\otimes \ab+ \one \otimes \ab\,\ab \\[2mm]
      \aabb &   \aabb\otimes \one + \ab \otimes \ab+ \one\otimes \aabb \\[2mm]
      \ab\,\aabb &  \ab\,\aabb \otimes \one + 2\,\ab\,\ab\otimes \ab + \ab\otimes \aabb + \ab \otimes \ab\,\ab+ \one \otimes \ab\,\aabb \\[2mm]
      \aabb\,\ab &\aabb\,\ab\otimes\one + \aabb\otimes\ab + \ab\otimes\ab\,\ab+ \one \otimes \aabb\,\ab  \\[2mm]
    \end{array}
  \end{equation*}
  \caption{Examples of the coproduct $\DeltaMKW$}\label{MKWcoprod}
\end{table}

The main result linking $\HMKW$ to LB-series is the following, which is an analog of the Hairer-Wanner theorem (Theorem \ref{hairerwannertheorem}) for B-series:
\begin{theorem}[{\cite{munthe-kaas2008oth}}]\label{compositionOfLB} The composition of two LB-series is again a LB-series: 
\begin{equation}
\Bs_f(\alpha)[\Bs_f(\beta)] = \Bs_f(\alpha * \beta),
\end{equation}
where $*$ is the convolution product in $\HMKW$.
\end{theorem}

\subsection{Lie--Butcher series and flows on manifolds.}
We shall see how LB-series can be used to represent numerical flows. More details and examples can be found in \cite{munthe-kaas1995lbt, munthe-kaas1998rkm, owren1999rkm, owren2006ocf, munthe-kaas2008oth, lundervold2009hao, lundervold2010bea}.

Flows $y_0 \mapsto y(t) = \Psi_t(y_0)$ on the manifold $\M$ can be
represented in several different ways. Here are three procedures, giving rise to what can be called LB-series of Type 1, 2 and 3:
\begin{enumerate}
\item In terms of pullback series: Find $\alpha\in G(\HMKW)$ such that 
  \begin{equation}\label{LBpullback}
\Psi(y(t)) = \Bs_t(\alpha)(y_0)[\Psi] \quad\mbox{for any $\Psi\in U(\g)^\M$.}
\end{equation}
This representation is used in the analysis of Crouch--Grossman methods by Owren and Marthinsen~\cite{owren1999rkm}.
In the classical setting, this is called a $S$-series~\cite{murua1999fsa}.
\item In terms of an autonomous differential equation: Find $\beta\in \g(\HMKW)$ such that $y(t)$ solves
\begin{equation}
y'(t) = \Bs_h(\beta)(y(t)).
\end{equation}
This is called backward error analysis (confer Section \ref{LGI:substitution}).
\item In terms of a non-autonomous equation of \emph{Lie type} (time dependent frozen vector field): Find $\gamma\in \g(\Hsh)$ such that $y(t)$ solves
\begin{equation}\label{lietype}
y'(t) = \left(\frac{\partial}{\partial t} \Bs_t(\gamma)(y_0)\right) y(t).
\end{equation}
This representation is used in~\cite{munthe-kaas1995lbt,munthe-kaas1998rkm}. In the classical setting this is (close to) the standard definition of $B$-series. 
\end{enumerate}
The algebraic relationships between the coefficients $\alpha$, $\beta$
and $\gamma$ in the above LB-series are \cite{lundervold2009hao}:
\begin{align*}
\beta &= \alpha\opr e &\mbox{$e$ is eulerian idempotent in $\HMKW$.}\\
\alpha &= \exp^\bpr(\beta)&\mbox{Exponential wrt.\ GL-product}\\
\gamma &= \alpha\opr Y^{-1}\opr D &\mbox{Dynkin idempotent in $\Hsh(\OT)$.}\\
\alpha &= Q(\gamma)&\mbox{$Q$-operator in $\Hsh(\OT)$.}
\end{align*}
The eulerian idempotent $e$ in a commutative, connected and graded
Hopf algebra $H$ is the formal series $e := \log^{*}(Id)$, where $Id$
is the identity endomorphism and $*$ the convolution product in $H$. The Dynkin map $D$ is the
convolution of the antipode $S$ and the grading operator $Y$, $D = S *
Y$, and $Y^{-1} \circ D$ is an idempotent. See e.g. \cite{lundervold2009hao} for details. The operator $Q$ is a rescaling of the Bell polynomials: 
\begin{align}
\begin{split}
Q_{n,k}(d_1,\ldots,d_{n-k+1}) &= \frac{1}{n!}B_{n,k}(1!d_1,\ldots,j!d_j,\ldots) = \mathop{\sum}_{|\omega| = n, \#(\omega)=k} \kappa({\omega})
\omega\\
Q_n(d_1,\ldots,d_n)&=\sum_{k=1}^n Q_{n,k}(d_1,\ldots,d_{n-k+1})\\
Q_0 &:= \one,
\end{split}
\end{align}
where, for $\omega = d_{j_1}d_{j_2}\cdots d_{j_k}$, the coefficients $\kappa({\omega})$ are defined as 
\[\kappa(\omega)  := \kappa(|d_{j_1}|,|d_{j_2}|,\ldots,|d_{j_k}|) := \frac{j_1j_2\cdots j_k}{j_1(j_1+j_2)\cdots(j_1+j_2+\cdots+j_k)}.\]

\noindent By using these relationships one can convert between the various
representations of flows.

\begin{example}[The exact solution]
The exact solution of a differential equation 
\[y'(t) = F(y(t))\] 
can be written as the solution of
\[y' = F_t\dpr y, \quad y(0)= y_0,  \]
where $F_t=F(y(t))\in \g$ is the pullback of $F$ along the time dependent flow of $F$. Let $F_t = \frac{\partial}{\partial t}\Bs_t(\gamma)$. By \cite[Proposition 4.9]{lundervold2009hao} the pullback is given by $\Bs_t(Q(\gamma_{\text{Exact}}))[F]$, so 
\[Y\opr\gamma_{\text{Exact}} = Q(\gamma_{\text{Exact}})[\ab]\Rightarrow \gamma_{\text{Exact}} = Y^{-1}\opr \Bplus(Q(\gamma_{\text{Exact}})).\]
Note that this is reminiscent of a so-called combinatorial Dyson--Schwinger equation~\cite{foissy2008fdb}. Solving by iteration yields
{\small
\begin{align*}
\gamma_{\text{Exact}} & = 
\ab + \frac{1}{2!}\aabb + \frac{1}{3!}\left(\aababb+\aaabbb\right) + \frac{1}{4!}\left(\aabababb+\aaabbabb+2\aabaabbb
+\aaababbb+\aaaabbbb\right)+ \frac{1}{5!}(\aababababb \\
&  +\aaabbababb+2\aabaabbabb +3\aababaabbb+\aaababbabb+ \aaaabbbabb+3\aaabbaabbb+3\aabaababbb+3\aabaaabbbb+\aaabababbb+\aaaabbabbb\\
& +2\aaabaabbbb+\aaaababbbb+\aaaaabbbbb)+\frac{1}{6!}\left(\aabababababb+\cdots\right)+\cdots
\end{align*}
}
\noindent A formula for the LB-series for the exact solution was given in \cite{owren1999rkm}. We observe that there cannot be any commutators of trees in this expression. Therefore, in LB-series of numerical integrators, commutators of trees must be zero up to the order of the method.
\end{example}

\begin{example}[The exponential Euler method]\label{eulerMethod}
The exponential Euler method \cite{iserles2000lgm} can be written as follows:
\[y_{n+1} = \exp(hf(y_n))y_n,\] or, by rescaling the vector field $f$, as \[y_{n+1} = \exp(f(y_n))y_n.\] This equation can be interpreted as a pullback equation of the form $\Phi(y_{n+1}) = \Bs(\exp(\ab))[\Phi]y_n$, so 
\[\alpha = \exp(\ab) = \one + \ab + \frac{1}{2!}\ab\ab + \frac{1}{3!}\ab\ab\ab + \cdots.\]
(Here the Grossman-Larson product is the same as concatenation). Note that $\exp(\ab) = Q(\ab),$ so the Type 3 LB-series for the Euler method is simply \[\gamma_{\text{Euler}} = \ab.\]
\end{example}

\begin{example}[The Lie--implicit midpoint method]\label{midpointMethod}
The Lie--implicit midpoint method \cite{iserles2000lgm} can be presented as:
\begin{align}\label{midpoint}
\sigma &= f(\exp(\frac{1}{2}\sigma)y_n) \\
y_{n+1} &= \exp(\sigma)y_n
\end{align}
We make the following ansatz: 
\begin{equation}\label{mp:ansatz}
\sigma = \sum_{\omega} \alpha(\omega)\omega = \alpha(\ab)\ab + \alpha\left(\aabb\right)\aabb + \alpha\left(\aababb\right)\aababb + \alpha\left([\ab,\aabb]\right)[\ab,\aabb] + \alpha\left(\aaabbb\right)\aaabbb + \cdots,
\end{equation}
i.e. that $\sigma$ can be written as an infinitesimal LB-series. From Equation (\ref{midpoint}), we get that
\begin{equation}\label{mp:sigma}
\sigma = \sum_{j=0}^{\infty} \frac{(\sigma)^j}{2^jj!} [\ab].
\end{equation}
Since there are no forests in this expression, we must have $\alpha([\omega, \omega']) = 0$ for all $\omega, \omega' \in \OT$. If we write $\tau= \Bplus(\tau_1 \cdots \tau_j)$ then, by combining Equation (\ref{mp:sigma}) with the ansatz, we see that coefficients of the LB-series are given recursively as $\alpha(\ab) = \frac{1}{2}$, 
\begin{equation}
\alpha(\tau) = \frac{1}{2^jj!} \alpha(\tau_1) \cdots \alpha(\tau_j).
\end{equation}
%The LB-series for the midpoint method is therefore \[\sum_{j, \tau = \Bplus(\tau_1 \cdots \tau_j)} \frac{1}{2^jj!} \alpha(\tau_1) \cdots \alpha(\tau_j) \mathcal{F}(\tau).\] 
Therefore,
\begin{align*}
\alpha_{\text{Midpoint}} & = 
\ab + \frac{1}{2!}\aabb + \frac{1}{2}\left(\frac{1}{4}\aababb+\aaabbb\right) +\cdots.
\end{align*}
\end{example}

\subsection{Substitution and backward error analysis for Lie--Butcher series}\label{LGI:substitution}
In \cite{lundervold2010bea} the substitution law for LB-series methods was developed, culminating in a formula that can be used to calculate the modified vector field used in backward error analysis. 

\paragraph{The substitution law.} The basic idea is as for B-series (Section \ref{substB}): We consider substituting a LB-series into another LB-series, e.g. $\Bs_{\Bs_f(\beta)}(\alpha)$, and the questions are as before: is this a LB-series, and in that case, which one? The result is given in terms of the \emph{substitution law}, defined using the freeness of the D-algebra $\N = \RR\langle \OT\rangle$ (Theorem \ref{universalD-alg}):

\begin{multicols}{2}
\begin{definition}\label{substlaw}For any map  $\alpha: \C \rightarrow D(\N)$ Theorem \ref{universalD-alg} implies that there a unique D-algebra homomorphism $\alpha\ast:\N\rightarrow \N$ such that $\alpha(c) = \alpha\ast c$ for all $c\in \C$. This homomorphism is called $\alpha$-substitution.

\vspace{-0.3cm}
\begin{diagram}[labelstyle=\scriptstyle]
\C &\rInto& \N\\
\dTo^{\alpha} && \dTo_{\alpha\ast} \\
D(\N) &\rInto& \N
\end{diagram}
\end{definition}
\end{multicols}

\begin{theorem}[{\cite{lundervold2010bea}}]\label{substBseries}
The substitution law defined in Definition \ref{substlaw} corresponds to the substitution of LB-series in the sense that 
\begin{equation}
\Bs_{\Bs_f(\alpha)}(\beta) = \Bs_f(\alpha \ast \beta).
\end{equation}
\end{theorem}

\paragraph{Calculating the substitution law.}
To obtain a formula for the substitution law we consider the dual $\astar$ of $\alpha$-substitution:
\begin{equation}
\langle \alpha \ast \beta, \omega \rangle = \langle \beta, \astar(\omega)\rangle,
\end{equation}
called the \emph{substitution character}. The dual pairing $\langle \cdot, \cdot \rangle$ is the one induced by requiring that all forests in $\OF$ are orthogonal, and we may write $\langle \alpha, \omega \rangle = \alpha(\omega)$. The map $\astar$ is a character for the shuffle product \cite{lundervold2010bea}: $\astar(\omega_1\shuffle\omega_2) = \astar(\omega_1) \shuffle \astar(\omega_2)$.

The formula for the substitution law is based on the cutting of trees as in the coproduct $\DeltaMKW$. More specifically, it is based on the dual of grafting, called \emph{pruning}:
\begin{equation}
\mathcal{P}_{\nu}(\omega) = \sum_{c \in LAC(\omega)} \langle \nu, P^c(\omega)\rangle R^c(\omega).
\end{equation}
Here the sum is over the left admissible cuts, but as opposed to the cuts in the formula (\ref{MKWcoprod_cuts}) for $\DeltaMKW$, the full cut is not included. 

In \cite{lundervold2010bea} the following inductive formula for $\astar$ was obtained:
\begin{theorem}[{\cite{lundervold2010bea}}]\label{substlawformula}
We have 
\begin{equation}\label{eq:substlaw}
\astar (\omega) = \sum_{(\omega) \in \DeltaC} \sum_{c \in LAC(\omega_{(2)})} \astar(\omega_{(1)}) \Bplus\left(\astar(P^c(\omega_{(2)}))\right)\alpha(R^c(\omega_{(2)})),
\end{equation}
if $\omega \neq 1$ and $\astar(\one) = \one$. Here $\DeltaC$ denotes the deconcatenation coproduct.
\end{theorem}
\noindent By using the magmatic operation $\times$ on $\N$, this can also be written as a composition of operators: 
\begin{equation}
\astar = \mu \circ (\mu_{\times} \otimes I)\circ (\astar \otimes \astar \otimes a) \circ (I\otimes \DeltaMKW') \circ \DeltaC.
\end{equation}
Here $\DeltaMKW'$ denotes the coproduct in (\ref{MKWcoprod_cuts}) with the full cut removed, and $\mu$ denotes concatenation.

Some examples of the substitution character can be found in Table \ref{subst_character}. Many more examples and details can be found in \cite{lundervold2010bea}. 
\begin{table}[h! ]
  \centering
  \begin{equation*}
    \begin{array}{c@{\,\,}|@{\quad}l}
      \hline \\[-2mm]
      \omega & \astar(\omega)  \\[1mm]
      \hline \\[-2mm]
      \one & \one \\[1mm]
      \ab & \alpha(\ab)\ab \\[2mm]
      \ab\ab & \alpha(\ab)^2\ab\ab\\[2mm]
      \aabb &   \alpha(\aabb)\ab + \alpha(\ab)^2\aabb\\[2mm]
      \ab\aabb &  \alpha(\ab\aabb)\ab + \alpha(\ab)\alpha(\aabb) \ab\ab + \alpha(\ab)^3\ab\aabb\\[2mm]
      \aabb\ab & \alpha(\aabb\ab)\ab + \alpha(\ab)\alpha(\aabb)\ab\ab + \alpha(\ab)^3 \aabb\ab \\[2mm]
    \end{array}
  \end{equation*}
  \caption{Examples of the substitution character $\astar$}\label{subst_character}
\end{table}

\begin{remark} One would like the substitution law $\ast$ for
  LB-series to be a convolution product in a Hopf or bialgebra, analogous to the substitution of B-series (Theorem \ref{thm:substB}). One way to achieve this is by obtaining a concrete description of the operations in the post-Lie operad. In that case one can follow the procedure in \cite{calaque2009tih}, which, roughly, is the following: The post-Lie operad has a pre-Lie structure (general phenomenon for augmented operads), there is an associated Lie algebra structure, its universal enveloping algebra is a Hopf algebra, and its dual is the Hopf algebra for the substitution law. This is a project currently under investigation \cite{ebrahimi-fard2011otp}.
\end{remark}

\section*{Acknowledgements}
We would like to thank Martin Bordemann, Kurusch Ebrahimi-Fard, Dominique Manchon and Ander Murua for many valuable remarks during discussions on the topics of this paper. Thanks also to H\aa{}kon Marthinsen and to the anonymous referees for their helpful comments. 

\bibliography{../../../bibliography/ref_alex}

\end{document}